\documentclass[11pt]{article}
\usepackage{latexsym,amssymb}
\usepackage{amsmath}
\usepackage[mathscr]{eucal}
\usepackage{color}
\usepackage[abbrev]{amsrefs}

\newtheorem{Theorem}{\bf Theorem}[section]
\newtheorem{Lemma}{\bf Lemma}[section]
\newtheorem{Proposition}{\bf Proposition}[section]
\newtheorem{Corollary}{\bf Corollary}[section]
\newtheorem{Remark}{\bf Remark}[section]
\newtheorem{Example}{\bf Example}[section]
\newtheorem{Definition}{\bf Definition}[section]

\newenvironment{theorem}{\begin{Theorem}$\!\!\!$}{\end{Theorem}}
\newenvironment{lemma}{\begin{Lemma}$\!\!\!$}{\end{Lemma}}
\newenvironment{proposition}{\begin{Proposition}$\!\!\!$}{\end{Proposition}}
\newenvironment{corollary}{\begin{Corollary}$\!\!\!$}{\end{Corollary}}
\newenvironment{remark}{\begin{Remark}$\!\!\!$}{\end{Remark}}

\newenvironment{definition}{\begin{Definition}$\!\!\!$}{\end{Definition}}

\numberwithin{equation}{section}

\numberwithin{equation}{section}

\newcommand{\dee}{{\rm{d}}}

\def\Xint#1{\mathchoice
{\XXint\displaystyle\textstyle{#1}}%
{\XXint\textstyle\scriptstyle{#1}}%
{\XXint\scriptstyle\scriptscriptstyle{#1}}%
{\XXint\scriptscriptstyle\scriptscriptstyle{#1}}%
\!\int}
\def\XXint#1#2#3{{\setbox0=\hbox{$#1{#2#3}{\int}$}
\vcenter{\hbox{$#2#3$}}\kern-.5\wd0}}

\def\dashint{\Xint-}

\allowdisplaybreaks

\begin{document}

\title{Initial traces and solvability of porous medium equation\\ with power nonlinearity}
\author{
\qquad\\
Kazuhiro Ishige, Nobuhito Miyake, and Ryuichi Sato
}
\date{}
\maketitle
\begin{abstract}
In this paper we study qualitative properties of initial traces of solutions to the porous medium equation 
with power nonlinearity, and obtain necessary conditions for the existence of solutions to the corresponding Cauchy problem. 
Furthermore, we establish sharp sufficient conditions
for the existence of solutions to the Cauchy problem using uniformly local Morrey spaces and their variations, 
and identify the optimal singularities of the initial data for the solvability of the Cauchy problem. 
\end{abstract}
\vspace{25pt}
\noindent Addresses:
\vspace{7pt}
\smallskip

\begin{tabular}{rl}
K.~I.: 
&Graduate School of Mathematical Sciences, The University of Tokyo,\\
&3-8-1 Komaba, Meguro-ku, Tokyo 153-8914, Japan.\\
&E-mail: {\tt ishige@ms.u-tokyo.ac.jp}\\[10pt]
N.~M.:  
&Faculty of Mathematics, Kyushu University, \\
&744, Motooka, Nishi-ku, Fukuoka 819-0395, Japan.\\
&E-mail: {\tt miyake@math.kyushu-u.ac.jp}\\[10pt]
R.~S.: 
& Faculty of Science  and Engineering, Iwate University,\\
& 4-3-5 Ueda, Morioka-shi, Iwate 020-8551, Japan.\\
& E-mail: {\tt ryusato@iwate-u.ac.jp}
\end{tabular}

\vspace{25pt}
\vspace{3pt}
%
\newpage
\section{Introduction}
Let $u$ be a nonnegative solution to the porous medium equation 
with power nonlinearity
\begin{equation}
\tag{E}
\label{eq:E}
\left\{
\begin{array}{ll}
\partial_t u=\Delta u^m+u^p & \quad\mbox{in}\quad{\mathbb R}^N\times(0,T),\vspace{3pt}\\
u\ge 0 & \quad\mbox{in}\quad{\mathbb R}^N\times(0,T),
\end{array}
\right.
\end{equation}
where $N\ge 1$, $1\le m<p$, $T\in(0,\infty]$, and $\partial_t:=\partial/\partial t$.
In this paper we study qualitative properties of initial traces of the solution~$u$, 
and obtain necessary conditions 
for the existence of solution to the Cauchy problem
\begin{equation}
\tag{P}
\label{eq:P}
\left\{
\begin{array}{ll}
\partial_t u=\Delta u^m+u^p & \quad\mbox{in}\quad{\mathbb R}^N\times(0,T),\vspace{3pt}\\
u(\cdot,0)=\mu & \quad\mbox{in}\quad{\mathbb R}^N,
\end{array}
\right.
\end{equation}
where $\mu$ is a nonnegative Radon measure in ${\mathbb R}^N$. 
Furthermore, we establish sharp sufficient conditions for the existence of nonnegative solutions 
to problem~\eqref{eq:P}, and identify the optimal singularities of the initial data 
for the solvability of problem~\eqref{eq:P}. 

In all that follows 
we denote by ${\mathcal M}$ the set of nonnegative Radon measures in ${\mathbb R}^N$. 
We also denote by $\mathcal{L}$ the set of nonnegative locally integrable functions in ${\mathbb R}^N$. 
We often identify $\dee\mu=\mu(x)\,\dee x$ in ${\mathcal M}$ for $\mu\in\mathcal{L}$. 
For any measurable set $E$ in ${\mathbb R}^d$, where $d=1,2,\dots$, 
let ${\mathcal L}^d(E)$ be the $d$-dimensional Lebesuge measure of $E$. 
For any $f\in\mathcal{L}$, $z\in{\mathbb R}^N$, and $r>0$, set
$$
\dashint_{B(z,r)} f\,\dee x:=\frac{1}{{\mathcal L}^N(B(z,r))}\int_{B(z,r)} f\,\dee x,
$$
where $B(z,r):=\{x\in{\mathbb R}^N\,:\,|x-z|<r\}$. 
Set
$$
p_m:=m+\frac{2}{N},\qquad \theta:=\frac{p-m}{2(p-1)},\qquad
\theta':=\frac{1}{\theta}=\frac{2(p-1)}{p-m}.
$$

The study of initial traces of nonnegative solutions to the Cauchy problem for parabolic equations is a classical subject, 
and qualitative properties of initial traces have been studied for various parabolic equations. 
See e.g.,~\cites{A} for linear parabolic equations, 
\cites{AC, HP} for porous medium equations, 
\cites{DH, DH02} for parabolic $p$\,-Laplace equations, 
\cites{I, IJK, ZX} for doubly nonlinear parabolic equations,
\cite{BSV} for fractional diffusion equations,
\cite{AIS} for Finsler heat equations,
\cites{ADi, FHIL, Hisa, HI18, HI24, HIT02, FI01, IKO, TY} for parabolic equations with source nonlinearity (positive nonlinearity),
\cites{BCV, BD, MV, MV02, MV03} for parabolic equations with absorption nonlinearity (negative nonlinearity). 

Let us recall some results on initial traces of solutions to problem~\eqref{eq:E} with $m=1$. 
See e.g.,  \cites{BP, HI18, IKO}. 
\begin{itemize}
  \item[(A)] 
  \begin{itemize}
  \item[(1)]
  Assume that problem~\eqref{eq:E} with $m=1$ possesses a solution~$u$ in ${\mathbb R}^N\times(0,T)$ for some $T\in(0,\infty)$. 
  Then there exists a unique $\nu\in{\mathcal M}$ such that 
  $$
  \underset{t\to+0}{\mbox{{\rm ess lim}}}\int_{{\mathbb R}^N} u(t)\psi\,\dee y=\int_{{\mathbb R}^N}\psi\,\dee \nu(y),\quad \psi\in C_c({\mathbb R}^N).
  $$
  Furthermore,  
  there exists $C_1=C_1(N,p)>0$ such that
  \begin{equation*}
  \sup_{z\in{\mathbb R}^N}\nu(B(z,\sigma))\le
  \left\{
  \begin{array}{ll}
  C_1 \sigma^{N-\frac{2}{p-1}} & \mbox{if}\quad p\not=p_1,\vspace{5pt}\\
  C_1\displaystyle{\left[\log\left(e+\frac{\sqrt{T}}{\sigma}\right)\right]^{-\frac{N}{2}}} & \mbox{if}\quad p=p_1,
  \end{array}
  \right.
  \end{equation*}
  for $\sigma\in(0,\sqrt{T})$. 
  \item[(2)] 
  Let $u$ be a solution to problem~\eqref{eq:P} with $m=1$.  
  Then $u$ is a solution to problem~\eqref{eq:E} and 
  the initial trace of $u$ coincides with the initial data of $u$ in ${\mathcal M}$.
  \end{itemize}
\end{itemize}
We remark that assertion~(A) gives necessary conditions for the existence of solutions to problem~\eqref{eq:P} with $m=1$. 
Sufficient conditions for the existence of solutions to problem~\eqref{eq:P} with $m=1$ 
have been studied in many papers, and the following assertion holds for $m=1$. 
\begin{itemize}
  \item[(B)] 
  Let $m=1$. 
  \begin{itemize}
  \item[(1)] Let $1<p<p_1$. 
  Then problem~\eqref{eq:P} possesses a local-in-time solution if and only if 
  $$
  \sup_{z\in{\mathbb R}^N}\mu(B(z,1))<\infty.
  $$
  (See assertion~(A) and e.g., \cites{W1, W2}.)
 \item[(2)] Let $p=p_1$.
  For any $\alpha>0$, 
  there exists $\epsilon_1=\epsilon_1(N,\alpha)>0$ such that 
  if $\mu\in\mathcal{L}$ satisfies
  $$
  \sup_{z\in{\mathbb R}^N}
  \sup_{\sigma\in(0,\sqrt{T})}
  \eta\left(\frac{\sigma}{\sqrt{T}}\right)
  \Psi^{-1}\left(\,\dashint_{B(z,\sigma)}
  \Psi\left(T^\frac{1}{p-1}\mu\right)\,\dee y\,\right)\le\epsilon_1
  $$
  for some $T\in(0,\infty)$, 
  then problem~\eqref{eq:P} possesses a solution in ${\mathbb R}^n\times(0,T)$,
  where
  $$
  \Psi(s):=s[\log (e+s)]^\alpha,
  \qquad
  \eta(s):=
  s^N\biggr[\log\biggr(e+\frac{1}{s}\biggr)\biggr]^{\frac{N}{2}}. 
  $$
  See e.g., \cites{FI01, HI18, IKO}. (See also \cite{IIK} for another sufficient condition.) 
  \item[(3)] Let $p>p_1$. 
  For any $\beta>1$, 
  there exists $\epsilon_2=\epsilon_2(N,p,\beta)>0$ such that
  if $\mu\in\mathcal{L}$ satisfies
  $$
  \sup_{z\in{\mathbb R}^n}\sup_{\sigma\in(0,\sqrt{T})}\,\sigma^{\frac{2}{p-1}}
  \left(\dashint_{B(z,\sigma)} |\mu|^\beta\,dy\,\right)^{\frac{1}{\beta}}\le\epsilon_2
  $$
  for some $T\in(0,\infty]$, then problem~\eqref{eq:P} possesses a solution in ${\mathbb R}^n\times(0,T)$. 
  See e.g., \cites{FI01, KY, HI18, RS}. 
   \end{itemize}
  \end{itemize}
We remark that assertion~(B)-(2) with $\alpha=0$ and assertion~(B)-(3) with $\beta=1$ do not hold.
(See \cite{T}*{Theorem~1}.) 
Combining the results in (A) and (B), we have 
the following assertion.
\begin{itemize}
\item[(C)]
Let $p\ge p_1$. 
For any $c>0$, set 
$$
\mu_c(x):=
\left\{
\begin{array}{ll}
c|x|^{-N}\displaystyle{\biggr[\log\biggr(e+\frac{1}{|x|}\biggr)\biggr]^{-\frac{N}{2}-1}} & \mbox{if}\quad p=p_1,\vspace{7pt}\\
c|x|^{-\frac{2}{p-1}} & \mbox{if}\quad p>p_1,
\end{array}
\right.
$$
for a.a.~$x\in{\mathbb R}^N$. 
\begin{itemize}
  \item[{\rm (1)}] 
  Problem~\eqref{eq:P} possesses a local-in-time solution for $c>0$ small enough; 
  \item[{\rm (2)}] 
  Problem~\eqref{eq:P} possesses no local-in-time solutions for $c>0$ large enough. 
\end{itemize}
Furthermore, if $p>p_1$ and $c>0$ is small enough, then 
problem~\eqref{eq:P} possesses a global-in-time solution. 
\end{itemize}
The results in (C) show that the ``strength" of the singularity at the origin of the functions $\mu_c$
is the critical threshold for the local solvability of problem~\eqref{eq:P}. 
We term such a singularity in the initial  data an {\em optimal singularity} of initial data
for the solvability of problem~\eqref{eq:P}. 
(See e.g., \cites{FHIL02} for further details of optimal singularities of  initial data.)
We easily see that, by translation invariance the singularity could be located at any point of ${\mathbb R}^N$. 

On the other hand, in the case of $m>1$, 
much less is known about the solvability of problem~\eqref{eq:P}, in particular, the optimal singularity 
of the initial data
for the solvability of problem~\eqref{eq:P}.
Andreucci and DiBenedetto~\cite{ADi} proved the existence and the uniqueness of 
initial traces of solutions to problem~\eqref{eq:E}. 
Furthermore, they obtained qualitative properties of initial traces of solutions and 
studied necessary conditions and sufficient conditions for the existence of solutions to problem~\eqref{eq:P}.
More precisely, they obtained the following assertion. 
\begin{itemize}
  \item[(D)] 
  Let $1\le m<p$.
  \begin{itemize}
  \item[(1)]
  Assume that problem~\eqref{eq:E} possesses a solution~$u$ in ${\mathbb R}^N\times(0,T)$ for some $T\in(0,\infty)$. 
  Then there exists a unique $\nu\in{\mathcal M}$ such that 
  $$
  \underset{t\to+0}{\mbox{{\rm ess lim}}}\int_{{\mathbb R}^N} u(t)\psi\,\dee y=\int_{{\mathbb R}^N}\psi\,\dee \nu(y),\quad \psi\in C_c({\mathbb R}^N).
  $$
  Furthermore, there exists $C_2=C_2(N,m,p)>0$ such that 
  \begin{equation*}
  \sup_{z\in{\mathbb R}^N}\nu(B(z,\sigma))\le C_2\sigma^{N-\frac{2}{p-m}}
  \end{equation*}
  for $\sigma\in(0,T^\theta)$. 
  \item[(2)] 
  Let $m<p<p_m$. 
  Then problem~\eqref{eq:P} possesses a local-in-time solution if and only if 
  $$
  \sup_{z\in{\mathbb R}^N}\mu(B(z,1))<\infty.
  $$
  \item[(3)] 
  Let $p\ge p_m$. 
  Then, for any $r>N(p-m)/2$, there exists $\epsilon_3=\epsilon_3(N,m, p,r)>0$ such that 
  if $\mu\in\mathcal{L}$ satisfies
  $$
  T^{\frac{1}{p-1}}
  \sup_{z\in{\mathbb R}^n}\left(\,\dashint_{B(z,T^\theta)}|\mu|^r\,\dee y\right)^{\frac{1}{r}}\le \epsilon_3
  $$
  for some $T\in(0,\infty)$, then problem~\eqref{eq:P} possesses a solution in ${\mathbb R}^N\times(0,T)$.
  \end{itemize}
\end{itemize}
Subsequently, the third author of this paper \cite{Sato} improved the result in assertion~(D)-(3) in the case of $p>p_m$, 
and obtained the following assertion.  
\begin{itemize}
  \item[(D)]
  \begin{itemize}
  \item[(3')]
  Let $p>p_m$. 
  Then there exists $\epsilon_4=\epsilon_4(N,m,p)>0$ such that 
  if $\mu\in\mathcal{L}$ satisfies
  $$
  T^{\frac{1}{p-1}}
  \sup_{z\in{\mathbb R}^n}\left(\,\dashint_{B(z,T^\theta)} |\mu|^{\frac{N(p-m)}{2}}\,\dee y\right)^{\frac{2}{N(p-m)}}\le \epsilon_4
  $$
  for some $T\in(0,\infty)$, then problem~\eqref{eq:P} possesses a solution in ${\mathbb R}^N\times(0,T)$. 
  \end{itemize}
\end{itemize}
Unfortunately, the results in (D) are not enough to identify
the optimal singularity of initial data for the solvability of problem~\eqref{eq:P} in the case of $p\ge p_m$.  
\vspace{3pt}

In this paper we give refinements of assertions~(D)-(1), (3), and (3'), 
and extend assertions~(A), (B), and (C) to the case $m>1$. 
More precisely: 
\begin{itemize}
	\item
		we improve qualitative properties of initial traces of solutions to problem~\eqref{eq:E} with $p=p_m$, 
		and establish sharp necessary condition for the existence of solutions to problem~\eqref{eq:P} 
		(see Theorem~\ref{Theorem:1.1});
	\item	
		we give sharp sufficient conditions for the existence of solutions to problem~\eqref{eq:P} with $p\ge p_m$ using uniformly local 
		Morrey spaces and their variations (see Theorems~\ref{Theorem:1.2} and \ref{Theorem:1.3}).
\end{itemize}
%
Our necessary conditions and sufficient conditions enable us to identify the optimal singularity of initial data for the solvability for problem~\eqref{eq:P} with $p\ge p_m$ (see Corollary~\ref{Corollary:1.1}).
\vspace{3pt}

We formulate definitions of solutions to problems~\eqref{eq:E} and \eqref{eq:P}. 
\begin{definition}
\label{Definition:1.1}
Let $1\le m<p$, $T\in(0,\infty]$, and $u\in L^p_{\rm loc}(\mathbb{R}^N\times[0, T))$ be nonnegative in ${\mathbb R}^N\times(0,T)$. 
\begin{itemize}
  \item[{\rm (1)}] 
  We say that $u$ is a solution to problem~\eqref{eq:E} in ${\mathbb R}^N\times(0,T)$ if $u$ satisfies 
  \begin{equation}
  \label{eq:1.1}
  \int^T_\tau\int_{\mathbb{R}^N}(-u\partial_t \phi-u^m\Delta\phi-u^p\phi)\,\dee x\,\dee t=\int_{\mathbb{R}^N}u(x,\tau)\phi(x,\tau)\,\dee x
  \end{equation}
  for $\phi\in C^{2;1}_c({\mathbb R}^N\times[0,T))$ and almost all {\rm({\it a.a.})} $\tau\in(0,T)$.
  \item[{\rm (2)}] 
  For any $\mu\in{\mathcal M}$,  
  we say that $u$ is a solution to problem~\eqref{eq:P} in ${\mathbb R}^N\times(0,T)$ if $u$ satisfies 
  $$
  \int^T_0 \int_{\mathbb{R}^N}(-u\partial_t \phi-u^m\Delta\phi-u^p\phi)\,\dee x\,\dee t=\int_{\mathbb{R}^N}\phi(x,0)\,\dee\mu(x)
  $$
  for $\phi\in C^{2;1}_c({\mathbb R}^N\times[0,T))$. 
\end{itemize}
\end{definition}

We introduce some notation. 
For any positive functions $f$ and $g$ in a set $X$, 
we say that $f\preceq g$ for $x\in X$ or equivalently that $g\succeq f$ for $x\in X$ if 
there exists $C>0$ such that 
$$
f(x)\le Cg(x)\quad\mbox{for $x\in X$}.
$$
If $f\preceq g$ and  $g\preceq f$ for $x\in X$, we say that $f\asymp g$ for $x\in X$.
In all that follows we will use $C$ to denote generic positive constants and point out that $C$  
may take different values  within a calculation. 

Let $\Phi$ be a nonnegative, convex, and strictly increasing function in $[0,\infty)$ such that $\Phi(0)=0$. 
Let $\rho$ be a nonnegative, non-decreasing, and continuous function in $[0,\infty)$. 
Then, for any $f\in\mathcal{L}$ and $R\in(0,\infty]$, 
set
$$
|||f|||_{\rho,\Phi;R}:=\sup_{z\in{\mathbb R}^N}\sup_{\sigma\in(0,R)}
\left\{\rho(\sigma)\Phi^{-1}\left(\,\dashint_{B(z,\sigma)}\Phi(f)\,\dee x\right)\right\}.
$$
For any $q\in[1,\infty)$ and $\alpha\in[1,\infty)$, 
if $\Phi(\xi)=\xi^\alpha$ and $\rho(\xi)=\xi^{\frac{N}{q}}$ for $\xi\in[0,\infty)$, 
we write 
$$
|||f|||_{q,\alpha;R}:=|||f|||_{\rho,\Phi;R}
=\sup_{z\in{\mathbb R}^N}\sup_{\sigma\in(0,R)}
\left\{\sigma^{\frac{N}{q}}\left(\,\dashint_{B(z,\sigma)}|f|^\alpha\,\dee x\right)^{\frac{1}{\alpha}}\right\}
$$
for simplicity. 
In particular, if $R=1$, then $|||\cdot|||_{q,\alpha;1}$ coincides with the norm of uniformly local Morrey spaces (see e.g., \cite{KY}*{Definition~0.1}).
\vspace{3pt}

Now we are ready to state the main results of this paper. 
The first theorem concerns with qualitative properties of initial traces of solutions to problem~\eqref{eq:E}. 
\begin{theorem}
\label{Theorem:1.1}
Let $N\ge 1$ and $1\le m<p$. 
\begin{itemize}
  \item[{\rm (1)}] 
  Let $u$ be a solution to problem~\eqref{eq:E} in ${\mathbb R}^N\times(0,T)$, where $T\in(0,\infty)$. 
  Then there exists a unique $\nu\in {\mathcal M}$ such that 
  \begin{equation}
  \label{eq:1.2}
  \underset{t\to+0}{\mbox{{\rm ess lim}}}\int_{{\mathbb R}^N}u(t)\psi\,\dee x=\int_{{\mathbb R}^N}\psi\,\dee\nu(x),
  \quad \psi\in C_c({\mathbb R}^N).
  \end{equation}
  Furthermore, there exists $C=C(N,m,p)>0$ such that
  \begin{equation}
  \label{eq:1.3}
  \sup_{z\in{\mathbb R}^N}\nu(B(z,\sigma))\le
  \left\{
  \begin{array}{ll}
  C\sigma^{N-\frac{2}{p-m}} & \mbox{if}\quad p\not=p_m,\vspace{5pt}\\
  \displaystyle{C\left[\log\left(e+\frac{T^{\theta}}{\sigma}\right)\right]^{-\frac{N}{2}}} & \mbox{if}\quad p=p_m,
  \end{array}
  \right.
  \end{equation}
  for $\sigma\in(0,T^{\theta})$.
  \item[{\rm (2)}] 
  Let $u$ be a solution to problem~\eqref{eq:P}. 
  Then $u$ is a solution to problem~\eqref{eq:E} and 
  the initial trace of $u$ coincides with the initial data of $u$ in ${\mathcal M}$.
\end{itemize}
\end{theorem}
Similarly to assertion~(A), Theorem~\ref{Theorem:1.1} gives a necessary condition for the existence of solutions to problem~\eqref{eq:P}. 
\begin{remark}
\label{Remark:1.1}
{\rm (1)} Let $u$ be a solution to problem~\eqref{eq:E} in ${\mathbb R}^N\times(0,T)$, where $T\in(0,\infty]$. 
For any $\lambda>0$, set 
$$
u_\lambda(x,t):=\lambda^{\frac{2}{p-m}}u(\lambda x,\lambda^{\theta'} t),\quad (x,t)\in{\mathbb R}^N\times(0,T_\lambda),
$$
where $T_\lambda:=\lambda^{-\theta'}T$. Then $u_\lambda$ is a solution to problem~\eqref{eq:E} in ${\mathbb R}^N\times(0,T_\lambda)$.
\vspace{3pt}
\newline
{\rm (2)} Let $1\le m<p<p_m$ and $\mu\in{\mathcal M}$. Then, since $N-2/(p-m)<0$, 
the relation 
$$
\sup_{z\in{\mathbb R}^N}\mu(B(z,\sigma))\le C\sigma^{N-\frac{2}{p-m}}
$$
holds for $\sigma\in(0,T^{\theta})$ if and only if 
$$
\sup_{z\in{\mathbb R}^N}\mu(B(z,T^{\theta}))\le CT^{N\theta-\frac{1}{p-1}}.
$$
{\rm (3)} 
If $m<p\le p_m$, then problem~\eqref{eq:P} possesses no nontrivial global-in-time solutions 
{\rm ({\it see} \cites{Gala, GKMS, EGKP, MP, MS, SGKM, Suzuki})}. 
This fact immediately follows from Theorem~{\rm\ref{Theorem:1.1}}. 
Indeed, assume that problem~\eqref{eq:P} possesses a nontrivial global-in-time solution~$u$. 
Then, for a.a.~$\tau>0$, setting $u_\tau(x,t):=u(x,t+\tau)$ for a.a.~$(x,t)\in{\mathbb R}^N\times(0,\infty)$, 
we see that $u_\tau$ is a global-in-time solution to problem~\eqref{eq:P} with initial data $\mu=u(\tau)$.
Combining this fact with Theorem~{\rm\ref{Theorem:1.1}}, we have
\begin{align*}
 & \sup_{z\in{\mathbb R}^N}\int_{B(z,T^{\theta})} u(\tau)\,\dee y\le CT^{N\theta-\frac{1}{p-1}}\to 0\quad\mbox{if}\quad p<p_m,\\
 & \sup_{z\in{\mathbb R}^N}\int_{B(z,T^{\frac{\theta}{2}})} u(\tau)\,\dee y
\le C\left[\log\left(e+T^{\frac{\theta}{2}}\right)\right]^{-\frac{N}{2}}\to 0\quad\mbox{if}\quad p=p_m,
\end{align*}
as $T\to\infty$. These imply that $u(x,\tau)=0$ for a.a.~$(x,\tau)\in{\mathbb R}^N\times(0,\infty)$, 
which is a contradiction. Thus problem~\eqref{eq:P} possesses no nontrivial global-in-time solutions if $m<p\le p_m$.
\end{remark}
In the second and third theorems
we obtain sharp sufficient conditions for the existence of solutions to problem~\eqref{eq:P} with $p=p_m$ 
and $p>p_m$, respectively. 
We remark that problem~\eqref{eq:P} possesses a global-in-time solution for some initial data $\mu\in\mathcal{L}$ 
if and only if $p>p_m$ (see Remark~\ref{Remark:1.1}-(3)).
\begin{theorem}
\label{Theorem:1.2}
Let $N\ge 1$, $m\ge 1$, $p=p_m$, and $\alpha>0$. 
Let
\begin{equation}
\label{eq:1.4}
\Psi(\xi):=\xi[\log (e+\xi)]^\alpha,
\qquad
\eta(\xi):=\xi^N\biggr[\log\biggr(e+\frac{1}{\xi}\biggr)\biggr]^{\frac{N}{2}},
\end{equation}
for $\xi\in[0,\infty)$. 
Then there exist $\epsilon_5=\epsilon_5(N,m,\alpha)>0$ and $C=C(N, m, \alpha)>0$ such that
if $\mu\in\mathcal{L}$ satisfies
\begin{equation}
\label{eq:1.5}
\sup_{z\in{\mathbb R}^N}\,\sup_{\sigma\in(0,T^{\theta})}\,
\eta\left(\frac{\sigma}{T^{\theta}}\right)
\Psi^{-1}\left(\,\dashint_{B(z,\sigma)}\Psi\left(T^{\frac{1}{p-1}}\mu\right)\,\dee y\right)
\le\epsilon_5
\end{equation}
for some $T\in(0,\infty)$,
then problem~\eqref{eq:P} possesses a solution $u$ in ${\mathbb R}^N\times(0,T)$, with $u$ satisfying
\begin{align*}
 & \sup_{t\in(0,T)}t^{\frac{1}{p-1}}\left[\log\left(e+\frac{T}{t}\right)\right]^{\frac{1}{p-1}}\|u(t)\|_{L^\infty({\mathbb R}^N)}\\
 & +\sup_{t\in(0,T)}\sup_{z\in{\mathbb R}^N}\,\sup_{\sigma\in(0,T^{\theta})}\,
\eta\left(\frac{\sigma}{T^{\theta}}\right)
\Psi^{-1}\left(\,\dashint_{B(z,\sigma)}\Psi\left(T^{\frac{1}{p-1}}u(t)\right)\,\dee y\right)\\
 & \le C\sup_{z\in{\mathbb R}^N}\,\sup_{\sigma\in(0,T^{\theta})}\,
\left\{\eta\left(\frac{\sigma}{T^{\theta}}\right)
\Psi^{-1}\left(\,\dashint_{B(z,\sigma)}\Psi\left(T^{\frac{1}{p-1}}\mu\right)\,\dee y\right)\right\}^{\frac{2}{N(m-1)+2}}.
\end{align*}
\end{theorem}
\begin{theorem}
\label{Theorem:1.3}
Let $N\ge 1$, $m\ge 1$, $p>p_m$, and $1<\beta<N(p-m)/2$. 
Then there exist $\epsilon_6=\epsilon_6(N,m,\beta)>0$ and $C=C(N, m, \beta)>0$ such that
if $\mu\in\mathcal{L}$ satisfies
\begin{equation}
\label{eq:1.6}
|||\mu|||_{\frac{N(p-m)}{2},\beta;\,T^{\theta}}\le\epsilon_6
\end{equation}
for some $T\in(0,\infty]$, 
then problem~\eqref{eq:P} possesses a solution $u$ in ${\mathbb R}^N\times(0,T)$, with $u$ satisfying
$$
\sup_{t\in(0,T)}t^{\frac{1}{p-1}}\|u(t)\|_{L^\infty}
+\sup_{t\in(0,T)}|||u(t)|||_{\frac{N(p-m)}{2},\beta;\,T^{\theta}}\le C|||\mu|||_{\frac{N(p-m)}{2},\beta;\,T^{\theta}}^{\frac{2\beta}{N(m-1)+2\beta}}.
$$
In particular, if $\mu\in\mathcal{L}$ satisfies
$$
|||\mu|||_{\frac{N(p-m)}{2},\beta;\,\infty}=\sup_{z\in{\mathbb R}^N}\sup_{\sigma>0}
\left\{\sigma^{\frac{2}{p-m}}\left(\,\dashint_{B(z,\sigma)}|\mu|^\beta\,\dee y\right)^{\frac{1}{\beta}}\right\}\le\epsilon_6,
$$
then problem~\eqref{eq:P} possesses a global-in-time solution. 
\end{theorem}

Necessary conditions in Theorem~\ref{Theorem:1.1} 
and sufficient conditions in Theorems~\ref{Theorem:1.2} and \ref{Theorem:1.3} are sharp. 
Indeed, Theorems~\ref{Theorem:1.1}--\ref{Theorem:1.3} enable us 
to identify an optimal singularity of the initial data for the solvability of problem~\eqref{eq:P} with $p\ge p_m$. 
\begin{corollary}
\label{Corollary:1.1}
Let $N\ge 1$, $m\ge 1$, and $p\ge p_m$. 
For any $c>0$, set 
$$
\mu_c^m(x):=
\left\{
\begin{array}{ll}
c|x|^{-N}\displaystyle{\biggr[\log\biggr(e+\frac{1}{|x|}\biggr)\biggr]^{-\frac{N}{2}-1}} & \mbox{if}\quad \displaystyle{p=p_m},\vspace{7pt}\\
c|x|^{-\frac{2}{p-m}} & \mbox{if}\quad \displaystyle{p>p_m},\vspace{3pt}
\end{array}
\right.
$$
for a.a.~$x\in{\mathbb R}^N$. 
\begin{itemize}
  \item[{\rm (1)}] 
  Problem~\eqref{eq:P} possesses a local-in-time solution for $c>0$ small enough; 
  \item[{\rm (2)}] 
  Problem~\eqref{eq:P} possesses no local-in-time solutions for $c>0$ large enough. 
\end{itemize}
Furthermore, if $p>p_m$ and $c>0$ is small enough, then 
problem~\eqref{eq:P} possesses a global-in-time solution.
\end{corollary}

In the proof of Theorem~\ref{Theorem:1.1} 
we modify arguments in \cite{TY} to find suitable cut-off functions to obtain estimates of solvable initial data of problem~\eqref{eq:P}
(see Proposition~\ref{Proposition:2.1}). 
Here we require delicate choice of parameters of cut-off functions (see Steps 3 and 4 in the proof of Proposition~\ref{Proposition:2.1}). 
Then we follow arguments in \cite{HI18} to complete the proof of Theorem~\ref{Theorem:1.1}.

The proofs of Theorems~\ref{Theorem:1.2} and \ref{Theorem:1.3} are new even if $m=1$. 
Indeed, the proofs of assertion~(B)-(2) in \cites{HI18, IIK, IKO}
and assertion~(B)-(3) in \cites{KY, HI18, RS} rely on the representation formula of solutions via Duhamel's principle. 
However, this approach is not available to the proofs of Theorems~\ref{Theorem:1.2} and \ref{Theorem:1.3} 
due to the nonlinearity of the principal term $\Delta u^m$.

One of the main ingredients of this paper is to give new energy estimates involving Morrey norms and its variations (see Lemmata~\ref{Lemma:3.1} and \ref{Lemma:3.2}).
The proofs of these lemmata are inspired by the argument in \cite{L}, which proved improved Sobolev inequalities via weak-type estimates and pseudo-Poincar\'{e} inequalities.
These energy estimates together with $L^\infty$-estimates of solutions to problem~\eqref{eq:P} 
lead a priori estimates for classical solutions to problem~\eqref{eq:P}.
Combining these estimates with regularity theorems for solutions to problem~\eqref{eq:P}, 
we construct a solution to problem~\eqref{eq:P}
as the limit of classical solutions to problem~\eqref{eq:P}, 
where initial data are given by the lifting and the truncation of $\mu$ (see \eqref{eq:4.6}).

The rest of this paper is organized as follows. 
In Section~2 we modify the arguments in \cite{TY} to prove Theorem~\ref{Theorem:1.1}. 
In Section~3 we obtain energy estimates of solutions. 
In Sections~4 and 5 we prove Theorems~\ref{Theorem:1.2} and \ref{Theorem:1.3}, respectively. 
In Section~6 we apply Theorems~\ref{Theorem:1.1}--\ref{Theorem:1.3} to prove Corollary~\ref{Corollary:1.1}.

\vspace{8pt}

\noindent
{\bf Acknowledgment.}
The authors of this paper would like to thank Professor Ryo Takada for his useful comments.
K. I. was supported in part by JSPS KAKENHI Grant Number 19H05599. 
N. M. was supported in part by JSPS KAKENHI Grant Numbers 22KJ0719 and 24K16944. 
R. S. was supported in part by JSPS KAKENHI Grant Number 21KK0044. 
\section{Proof of Theorem~\ref{Theorem:1.1}}
%
In this section we modify arguments in \cites{TY} to study necessary conditions for the existence of solutions to problem~\eqref{eq:P}. 
Furthermore, we obtain qualitative properties of initial traces of solutions to problem~\eqref{eq:E}, 
and prove Theorem~\ref{Theorem:1.1}.
\begin{proposition}
\label{Proposition:2.1}
Let $u$ be a solution to problem~\eqref{eq:P} in ${\mathbb R}^N\times[0,T)$, where $T\in(0,\infty)$.  
Then there exists $C=C(N,m,p)>0$ such that
$$
\sup_{z\in{\mathbb R}^N}\mu(B(z,\sigma))\le
\left\{
\begin{array}{ll}
C\sigma^{N-\frac{2}{p-m}} & \mbox{if}\quad p\not=p_m,\vspace{5pt}\\
\displaystyle{C\left[\log\left(e+\frac{T^{\theta}}{\sigma}\right)\right]^{-\frac{N}{2}}} & \mbox{if}\quad p=p_m,
\end{array}
\right.
$$
for $\sigma\in(0,T^{\theta})$.
\end{proposition}
{\bf Proof.} 
The proof is divided into several steps. 
Let $u$ be a solution to problem~\eqref{eq:P} in ${\mathbb R}^N\times(0,T)$, where $T\in(0,\infty)$.
\vspace{3pt}
\newline
\underline{Step 1}: 
Let $\psi\in C_c^{2;1}({\mathbb R}^N\times[0,T))$ be chosen later such that $0\le\psi\le 1$ in ${\mathbb R}^N\times[0,T)$. 
Let $k\in{\mathbb N}$, and set $\phi=\psi^k$. 
Then it follows from Definition~\ref{Definition:1.1}-(2) and H\"older's inequality that
\begin{alignat*}{1}
	&\int_{\mathbb{R}^N}\phi(0)\,\dee\mu(x)+\int^T_0\int_{\mathbb{R}^N}u^p\phi\,\dee x\,\dee t
	=-\int^T_0\int_{\mathbb{R}^N}(u\partial_t \phi+u^m\Delta\phi)\,\dee x\,\dee t\\
	&\quad \le\bigg(\int^T_0\int_{\mathbb{R}^N}u^p\phi\,\dee x\,\dee t\bigg)^{\frac{1}{p}}\bigg(\int^T_0\int_{\mathbb{R}^N}
	\left|\frac{\partial_t\phi}{\phi}\right|^{\frac{p}{p-1}}\phi\,\dee x\,\dee t
	\bigg)^{1-\frac{1}{p}}\\
	&\qquad\qquad
	+\bigg(\int^T_0\int_{\mathbb{R}^N}u^p\phi\,\dee x\,\dee t\bigg)^{\frac{m}{p}}\bigg(\int^T_0\int_{\mathbb{R}^N}
	\left|\frac{\Delta\phi}{\phi}\right|^{\frac{p}{p-m}}\phi\,\dee x\,\dee t
	\bigg)^{1-\frac{m}{p}}\\
	&\le\int^T_0\int_{\mathbb{R}^N}u^p\phi\,\dee x\,\dee t
	+C\int^T_0\int_{\mathbb{R}^N}\left(|\partial_t \phi|^{\frac{p}{p-1}}\phi^{-\frac{1}{p-1}}+|\Delta\phi|^{\frac{p}{p-m}}\phi^{-\frac{m}{p-m}}\right)\,\dee x\,\dee t
\end{alignat*}
and hence
\begin{equation}
\label{eq:2.1}
	\int_{\mathbb{R}^N}\phi(0)\,\dee\mu(x)
	\le C\int^T_0\int_{\mathbb{R}^N}\left(|\partial_t \phi|^{\frac{p}{p-1}}\phi^{-\frac{1}{p-1}}+|\Delta\phi|^{\frac{p}{p-m}}\phi^{-\frac{m}{p-m}}\right)\,\dee x\,\dee t.
\end{equation}
Since $\psi\le 1$ in ${\mathbb R}^N\times[0,\infty)$ and $p>m\ge 1$,
taking $k\ge 1$ large enough, we see that
\begin{alignat*}{1}
	 & |\partial_t \phi|^{\frac{p}{p-1}}\phi^{-\frac{1}{p-1}}
	\le C|\partial_t \psi|^{\frac{p}{p-1}}\psi^{\frac{k(p-1)-p}{p-1}}\le C|\partial_t \psi|^{\frac{p}{p-1}},
	\\
	 & |\Delta\phi|^{\frac{p}{p-m}}\phi^{-\frac{m}{p-m}}
	\le C|\Delta\psi|^{\frac{p}{p-m}}\psi^{\frac{k(p-m)-p}{p-m}}+|\nabla\psi|^{\frac{2p}{p-m}}\psi^{\frac{k(p-m)-2p}{p-m}}
	\le C|\Delta\psi|^{\frac{p}{p-m}}+C|\nabla\psi|^{\frac{2p}{p-m}},
\end{alignat*}
which together with \eqref{eq:2.1} imply that
\begin{equation}
\label{eq:2.2}
\begin{split}
\int_{\{x\in\mathbb{R}^N\,:\,\psi(x, 0)=1\}}\,\dee\mu(x)
\le C\int^T_0\int_{\mathbb{R}^N}\Big(|\partial_t \psi|^{\frac{p}{p-1}}+|\Delta\psi|^{\frac{p}{p-m}}+|\nabla\psi|^{\frac{2p}{p-m}}\Big)\,\dee x\,\dee t\\
\le C\int^T_0\int_{\mathbb{R}^N}\Big(|\partial_t \psi|^{\frac{p}{p-1}}+(|\Delta\psi|+|\nabla\psi|^2)^{\frac{p}{p-m}}\Big)\,\dee x\,\dee t.
\end{split}
\end{equation}
\underline{Step 2}: 
Let $\zeta\in C^\infty(\mathbb{R})$ be such that 
$0\le \zeta\le 1$ in $\mathbb{R}$, $\zeta\equiv 1$ in $[1, \infty)$, $\zeta\equiv 0$ in $(-\infty, 0]$, 
and $|\zeta'|\le 2$ in $\mathbb{R}$. 
Let $z\in{\mathbb R}^N$, $a>0$, and $\delta>0$. 
For any $F\in C^\infty((0,\infty))$ with $F\le 0$ in $[aT/2,\infty)$, set 
$$
\psi(x,t):=\zeta(F(r(x,t)))\quad\mbox{with}\quad
r(x,t):=|x-z|^{\theta'}+at+\delta
$$
for $(x,t)\in {\mathbb R}^N\times[0,T)$. 
Then $\psi\in C_c^2({\mathbb R}^N\times[0,T))$. 
Since 
\begin{align*}
 & \theta'=\dfrac{2(p-1)}{p-m}>2,\\
 & \partial_t\psi=\zeta'(F(r))F'(r)\partial_t r,
\quad
\nabla\psi=\zeta'(F(r))F'(r)\nabla r,\\
 & \Delta\psi=\zeta''(F(r))F'(r)^2|\nabla r|^2
+\zeta'(F(r))F''(r)|\nabla r|^2
+\zeta'(F(r))F'(r)\Delta r,
\end{align*}
we have
\begin{align*}
|\partial_t\psi| & \le Ca|F'(r)|,\\
|\nabla\psi|^2 & \le CF'(r)^2|x-z|^{2\left(\theta'-1\right)}\le CF'(r)^2r^{2-2\theta},\\
|\Delta\psi| & \le CF'(r)^2|x-z|^{2\left(\theta'-1\right)}+C|F''(r)||x-z|^{2\left(\theta'-1\right)}+C|F'(r)||x-z|^{\theta'-2}\\
 & \le CF'(r)^2r^{2-2\theta}+C|F''(r)|r^{2-2\theta}+C|F'(r)|r^{1-2\theta}.
\end{align*}
These together with \eqref{eq:2.2} imply that
\begin{equation*}
\int_{\{x\in\mathbb{R}^N\,:\,F(|x-z|^{\theta'}+\delta)\ge 1\}}\,\dee\mu(x)\le C\iint_{\{(x,t)\in{\mathbb R}^N\times[0,T)\,:\,0\le F(r(x,t))\le 1\}} g(r(x,t))\,\dee x\,\dee t,
\end{equation*}
where
\begin{equation}
\label{eq:2.3}
g(\xi):=a^{\frac{p}{p-1}}|F'(\xi)|^{\frac{p}{p-1}}+\Big(|F'(\xi)|^2\xi^{2-2\theta}+|F''(\xi)|\xi^{2-2\theta}+|F'(\xi)|\xi^{1-2\theta}\Big)^{\frac{p}{p-m}}.
\end{equation}
Then we obtain
\begin{equation}
\label{eq:2.4}
\begin{split}
 & \sup_{z\in{\mathbb R}^N}\int_{\{x\in\mathbb{R}^N\,:\,F(|x-z|^{\theta'}+\delta)\ge 1\}}\,\dee\mu(x)\\
 & \le \sup_{z\in\mathbb{R}^N}\iint_{\{(x,t)\in{\mathbb R}^N\times[0,\infty)\,:\,0\le F(|x-z|^{\theta'}+at+\delta)\le 1\}} g(|x-z|^{\theta'}+at+\delta)\,\dee x\,\dee t\\
 & =Ca^{-1}\iint_{\{(r,t)\in[0,\infty)\times[0,\infty)\,:\,0\le F(r^{\theta'}+t+\delta)\le 1\}} g(r^{\theta'}+t+\delta)r^{N-1}\,\dee r\,\dee t\\
 & =Ca^{-1}\iint_{\{(s,\tau)\in[0,\infty)\times[0,\infty)\,:\,0\le F(s^2+\tau^2+\delta)\le 1\}} g(s^2+\tau^2+\delta)s^{2(N-1)\theta+2\theta-1}\tau\,\dee s\,\dee \tau\\
 & 
 =Ca^{-1}\int_{\{\zeta\ge 0\,:\, 0\le F(\zeta^2+\delta)\le 1\}}g(\zeta^2+\delta) \zeta^{2N\theta+1}\,\dee\zeta\,
 \int^\frac{\pi}{2}_0\cos^{2N\theta-1}(\omega)\sin(\omega)\,\dee \omega
 \\
 & =Ca^{-1}\int_{\{\xi\ge 0\,:\, 0\le F(\xi+\delta)\le 1\}}g(\xi+\delta) \xi^{N\theta}\,\dee\xi.
\end{split}
\end{equation}
\underline{Step 3}: 
Let $b$, $c$, and $d>0$ be constants to be chosen later such that
\begin{equation}
\label{eq:2.5}
\frac{aT}{2}\ge \dfrac{d}{e^c-1}.
\end{equation}
Set
\[
F(\xi):=\dfrac{1}{b}\left(\log\left(1+\dfrac{d}{\xi}\right)-c\right)\quad\mbox{for}\quad \xi\in(0,\infty).
\]
Then 
\begin{equation}
\label{eq:2.6}
\begin{split}
 & F(\xi)\ge 1\quad\mbox{if and only if}\quad \xi\le R_1:=\frac{d}{e^{b+c}-1},\\
 & F(\xi)\ge 0\quad\mbox{if and only if}\quad \xi\le R_2:=\frac{d}{e^c-1},\\
 & F(\xi)\le 0\quad\mbox{if}\quad \xi\ge aT/2. 
\end{split}
\end{equation}
Furthermore, 
\begin{equation}
\label{eq:2.7}
\begin{split}
 & |F'(\xi)|=\left|-\frac{1}{b}\left(1+\frac{d}{\xi}\right)^{-1}\frac{d}{\xi^2}\right|=\frac{d}{b\xi}\frac{1}{\xi+d}\le\frac{1}{b\xi},\\
 & |F''(\xi)|=\frac{d}{b\xi^2}\frac{1}{\xi+d}+\frac{d}{b\xi}\frac{1}{(\xi+d)^2}
=\frac{1}{b\xi^2}\frac{d(\xi+d)+d\xi}{(\xi+d)^2}\le\frac{1}{b\xi^2},
\end{split}
\end{equation}
for $\xi\in(0,\infty)$. 
Letting $\delta\to 0$ and applying \eqref{eq:2.4}, 
by \eqref{eq:2.3}, \eqref{eq:2.6}, and \eqref{eq:2.7} we obtain 
\begin{equation}
\label{eq:2.8}
\begin{split}
 & \sup_{z\in{\mathbb R}^N}\mu(B(z,R_1^{\theta}))\\
 & \le C_*a^{-1}\big(a^{\frac{p}{p-1}}b^{-\frac{p}{p-1}}+b^{-\frac{2p}{p-m}}+b^{-\frac{p}{p-m}}\big)\int^{R_2}_{R_1}\xi^{N\theta-\frac{p}{p-1}}\,\dee\xi\\
 & =C_*R_1^{N\theta-\frac{1}{p-1}} 
 a^{-1}\big(a^{\frac{p}{p-1}}b^{-\frac{p}{p-1}}+b^{-\frac{2p}{p-m}}+b^{-\frac{p}{p-m}}\big)\int^{R_2/R_1}_1\xi^{N\theta-\frac{p}{p-1}}\,\dee\xi,
 \end{split}
 \end{equation}
where $C_*$ is a positive constant independent of $a$, $b$, $c$, and $d$. 

Let $\sigma>0$.  
Let $d>0$ be such that 
$$
\sigma=\left(\dfrac{d}{e^{b+c}-1}\right)^{\theta}=R_1^{\theta}.
$$
Then it follows from \eqref{eq:2.5} and \eqref{eq:2.8} that
$$
\sup_{z\in{\mathbb R}^N}\mu(B(z,\sigma))
\le C_*\sigma^{N-\frac{2}{p-m}}a^{-1}\big(a^{\frac{p}{p-1}}b^{-\frac{p}{p-1}}+b^{-\frac{2p}{p-m}}+b^{-\frac{p}{p-m}}\big)\int^{\frac{e^{b+c}-1}{e^{c}-1}}_{1} \xi^{N\theta-\frac{p}{p-1}}\,\dee\xi
$$
for $a$, $b$, $c$, and $\sigma>0$ with
$$
0<\sigma\le \left(\dfrac{aT(e^c-1)}{2(e^{b+c}-1)}\right)^{\theta}.
$$
Letting $c\to\infty$, we get
\begin{equation}
\label{eq:2.9}
\sup_{z\in{\mathbb R}^N}\mu(B(z,\sigma))
\le C_*\sigma^{N-\frac{2}{p-m}}a^{-1}\big(a^{\frac{p}{p-1}}b^{-\frac{p}{p-1}}+b^{-\frac{2p}{p-m}}+b^{-\frac{p}{p-m}}\big)\int^{e^{b}}_{1} \xi^{N\theta-\frac{p}{p-1}}\,\dee\xi
\end{equation}
for $a$, $b$, and $\sigma>0$ with $0<\sigma\le (aT/2e^b)^{\theta}$. 
\vspace{3pt}
\newline
\underline{Step 4}: 
We choose suitable $a$ and $b>0$ to complete the proof of Proposition~\ref{Proposition:2.1}. 
In the case of $p\not=p_m$, 
by \eqref{eq:2.9} with $a=2e$ and $b=1$ 
we obtain 
$$
\mu(B(z,\sigma))
\le C\sigma^{N-\frac{2}{p-m}}\quad\mbox{for}\quad \sigma\in(0,T^{\theta}). 
$$
Thus Proposition~\ref{Proposition:2.1} follows in the case of $p\not=p_m$. 

Consider the case of $p=p_m$. 
Let $\ell\ge 1$ and $b\ge 1$, and set
$a=\ell b^{-\frac{m-1}{p-m}}$. 
Since
\begin{align*}
 & a^{\frac{p}{p-1}}b^{-\frac{p}{p-1}}=\left(\ell b^{-\frac{p-1}{p-m}}\right)^{\frac{p}{p-1}}=\ell^{\frac{p}{p-1}}b^{-\frac{p}{p-m}},\\
 & N\theta-\frac{p}{p-1}=-1, 
\quad
\frac{m-1}{p-m}-\frac{p}{p-m}+1=-\frac{1}{p-m}=-\frac{N}{2},
\end{align*}
it follows from \eqref{eq:2.9} that 
\begin{equation}
\label{eq:2.10}
\sup_{z\in{\mathbb R}^N}\mu(B(z,\sigma))\le C\ell^{\frac{p}{p-1}}a^{-1}b^{-\frac{p}{p-m}+1}=C\ell^{\frac{1}{p-1}}b^{-\frac{N}{2}}
\end{equation}
for $b\ge 1$ and $\sigma>0$ with 
$$
0<\sigma\le\Biggr(\,\frac{\ell T}{2b^{\frac{N(m-1)}{2}}e^b}\,\Biggr)^{\theta}.
$$ 
Let $L\ge e$, and set
$$
b=\log\left(\left(L+\frac{T}{\sigma^{\theta'}}\right)\left[\log\left(L+\frac{T}{\sigma^{\theta'}}\right)\right]^{-\frac{N(m-1)}{2}}\right).
$$
Taking $L$ large enough if necessary, we see that $b\ge 1$. 
Since
$$
b\le \log\left(L+\frac{T}{\sigma^{\theta'}}\right)\quad\mbox{for $\sigma>0$},
$$
we have
$$
b^{\frac{N(m-1)}{2}}e^b\le \left[\log\left(L+\frac{T}{\sigma^{\theta'}}\right)\right]^{\frac{N(m-1)}{2}}
\left(L+\frac{T}{\sigma^{\theta'}}\right)\left[\log\left(L+\frac{T}{\sigma^{\theta'}}\right)\right]^{-\frac{N(m-1)}{2}}
=L+\frac{T}{\sigma^{\theta'}}\le C\frac{T}{\sigma^{\theta'}}
$$
for $\sigma\in(0,T^{\theta})$. 
Then, taking $\ell\ge 1$ large enough if necessary, we obtain
$$
\frac{\ell T}{2b^{\frac{N(m-1)}{2}}e^b}\ge\frac{\ell}{C}\sigma^{\theta'}\ge \sigma^{\theta'}\quad\mbox{for}\quad \sigma\in(0,T^{\theta}),
$$
which implies that
$$
\Biggr(\,\frac{\ell T}{2b^{\frac{N(m-1)}{2}}e^b}\,\Biggr)^{\theta}\ge\sigma \quad\mbox{for}\quad \sigma\in(0,T^{\theta}). 
$$
Then we deduce from \eqref{eq:2.10} that
\begin{align*}
\sup_{z\in{\mathbb R}^N}\mu(B(z,\sigma)) & \le C\left[\log\left(\left(L+\frac{T}{\sigma^{\theta'}}\right)\left[\log\left(L+\frac{T}{\sigma^{\theta'}}\right)\right]^{-\frac{N(m-1)}{2}}\right)\right]^{-\frac{N}{2}}\\
 & \le C\left[\log\left(L+\frac{T}{\sigma^{\theta'}}\right)\right]^{-\frac{N}{2}}
\le C\left[\log\left(e+\frac{T^{\theta}}{\sigma}\right)\right]^{-\frac{N}{2}}
\end{align*}
for $\sigma\in(0,T^{\theta})$. 
Thus Proposition~\ref{Proposition:2.1} follows in the case of $p=p_m$, 
and the proof of Proposition~\ref{Proposition:2.1} is complete. 
$\Box$
\vspace{5pt}
\newline
{\bf Proof of Theorem~\ref{Theorem:1.1}.}
Let $u$ be a solution to problem~\eqref{eq:E} in ${\mathbb R}^N\times(0,T)$, where $T\in(0,\infty)$.
Then there exists a measurable set $I\subset(0,T)$ with ${\mathcal L}^1((0,T)\setminus I)=0$ such that 
\eqref{eq:1.1} holds for $\tau\in I$. 
For any $\tau\in I$, setting $u_\tau(x,t):=u(x,t+\tau)$ for a.a.~$(x,t)\in{\mathbb R}^N\times(0,T-\tau)$, 
we see that $u_\tau$ is a solution to problem~\eqref{eq:P} in ${\mathbb R}^N\times(0,T-\tau)$ with $\mu=u(\tau)$. 
Then, by Proposition~\ref{Proposition:2.1} we see that 
\begin{equation}
\label{eq:2.11}
\sup_{z\in{\mathbb R}^N}\int_{B(z,\sigma)}u(\tau)\,\dee x\le
\left\{
\begin{array}{ll}
C_*\sigma^{N-\frac{2}{p-m}} & \mbox{if}\quad p\not=p_m,\vspace{5pt}\\
\displaystyle{C_*\left[\log\left(e+\frac{T^{\theta}}{\sigma}\right)\right]^{-\frac{N}{2}}} & \mbox{if}\quad p=p_m,
\end{array}
\right.
\end{equation}
for $\sigma\in(0,(T-\tau)^{\theta})$ and $\tau\in I$, 
where $C_*$ is a positive constant depending only on $N$, $p$, and $m$. 
Applying the weak compactness of Radon measures (see e.g., \cite{EG}*{Section~1.9}),
we find a sequence $\{\tau_j\}\subset I$ with $\lim_{j\to\infty}\tau_j=0$ and $\nu\in{\mathcal M}$ such that 
\begin{equation}
\label{eq:2.12}
\lim_{j\to\infty}\int_{{\mathbb R}^N} u(\tau_j)\psi\,\dee x=\int_{{\mathbb R}^N}\psi\,\dee \nu(x),
\quad \psi\in C_c({\mathbb R}^N).
\end{equation}

We show that \eqref{eq:1.2} holds. 
Let $\{s_j\}\subset I$ with $\lim_{j\to\infty}s_j=0$ and $\nu'\in{\mathcal M}$ such that 
\begin{equation}
\label{eq:2.13}
\lim_{j\to\infty}\int_{{\mathbb R}^N} u(s_j)\psi\,\dee x=\int_{{\mathbb R}^N}\psi\,\dee \nu'(x),
\quad \psi\in C_c({\mathbb R}^N).
\end{equation}
Let $\psi\in C^\infty_c({\mathbb R}^N)$. 
Let $\phi\in C_c^\infty({\mathbb R}^N\times[0,T))$ be such that $\phi(x,t)=\psi(x)$ for $(x,t)\in{\mathbb R}^N\times[0,\delta]$ for some $\delta\in(0,T)$. 
Then, by \eqref{eq:1.1}, \eqref{eq:2.12}, and \eqref{eq:2.13} we see that 
\begin{align*}
\int^T_0 \int_{\mathbb{R}^N}(-u\partial_t \phi-u^m\Delta\phi-u^p\phi)\,\dee x\,\dee t
 & =\lim_{j\to\infty}\int_{{\mathbb R}^N}u(\tau_j)\phi(\tau_j)\,\dee x
=\int_{\mathbb{R}^N}\psi\,\dee\nu(x)\\
 & =\lim_{j\to\infty}\int_{{\mathbb R}^N}u(s_j)\phi(s_j)\,\dee x
=\int_{\mathbb{R}^N}\psi\,\dee\nu'(x).
\end{align*}
This implies that $\nu=\nu'$ in ${\mathcal M}$. 
Then, since $\{s_j\}\subset I$ is arbitrary,
we see that $\nu$ satisfies \eqref{eq:1.2}. Furthermore, the uniqueness of the initial trace of solution~$u$ also follows. 

It remains to prove \eqref{eq:1.3}. 
Let $z\in{\mathbb R}^N$, $\delta\in(0,T)$, $\sigma\in(0,(T-\delta)^{\theta})$, and $\epsilon\in(0,\sigma)$. 
Let $\psi\in C_c({\mathbb R}^N)$ be such that $0\le\psi\le 1$ in ${\mathbb R}^N$, $\psi=1$ in $B(z,\sigma-\epsilon)$, and $\mbox{supp}\,\psi\subset B(z,\sigma)$.
Then it follows from \eqref{eq:2.11} that
\begin{align*}
\mu(B(z,\sigma-\epsilon)) & \le\int_{{\mathbb R}^N}\psi\,\dee\mu(x)=\lim_{j\to 0}\int_{{\mathbb R}^N} u(\tau_j)\psi\,\dee x\\
 & \le
\left\{
\begin{array}{ll}
C_*\sigma^{N-\frac{2}{p-m}} & \mbox{if}\quad p\not=p_m,\vspace{5pt}\\
\displaystyle{C_*\left[\log\left(e+\frac{T^{\theta}}{\sigma}\right)\right]^{-\frac{N}{2}}} & \mbox{if}\quad p=p_m.
\end{array}
\right.
\end{align*}
Since $\epsilon\in(0,\sigma)$ and $\delta\in(0,T)$ are arbitrary, 
we obtain \eqref{eq:1.3}. 
Thus Theorem~\ref{Theorem:1.1} follows.
$\Box$
\section{Energy estimates of solutions}
In this section we obtain energy estimates of solutions to problem~\eqref{eq:P}, 
which are crucial in the proofs of Theorems~\ref{Theorem:1.2} and \ref{Theorem:1.3}. 
We often use the following property: 
\begin{itemize}
  \item there exists $m_*\ge1$ such that
  \begin{equation}
  \label{eq:3.1}
  \sup_{z\in{\mathbb R}^N}\int_{B(z,2\sigma)}f\,\dee x\le m_*\sup_{z\in{\mathbb R}^N}\int_{B(z,\sigma)}f\,\dee x
  \end{equation}
  for $f\in\mathcal{L}$ and $\sigma>0$ (see e.g., \cite{IS01}*{Lemma 2.1}). 
\end{itemize}
We first give an energy estimate of solutions to problem~\eqref{eq:P} with $p=p_m$. 
\begin{lemma}
\label{Lemma:3.1}
Let $p=p_m$, $\mu\in \mathcal{L}\cap L^\infty({\mathbb R}^N)$, and $T\in(0,\infty)$. 
Let $u$ be a positive classical solution to problem~\eqref{eq:P} in ${\mathbb R}^N\times(0,T)$ 
such that 
\begin{equation}
\label{eq:3.2}
\sup_{t\in(0,T)}\|u(t)\|_{L^\infty({\mathbb R}^N)}<\infty.
\end{equation}
Let $\Psi$ and $\eta$ be as in \eqref{eq:1.4}. 
Then 
\begin{equation}
\label{eq:3.3}
\sup_{z\in{\mathbb R}^N}\int^t_0\int_{B(z,\sigma)}u^{m-1}\Psi''(u)|\nabla u|^2\,\dee x\,\dee s<\infty
\end{equation}
for $t\in(0,T)$ and $\sigma>0$.
Furthermore, there exists $C=C(N,m,\alpha)>0$ such that
\begin{equation}
\label{eq:3.4}
\begin{split} 
 & \sup_{s\in(0,t]}\sup_{z\in{\mathbb R}^N}\int_{B(z,\sigma)}\Psi(u(s))\,\dee x+\sup_{z\in{\mathbb R}^N}\int^t_0\int_{B(z,\sigma)}u^{m-1}\Psi''(u)|\nabla u|^2\,\dee x\,\dee s\\
 & \le C\sup_{z\in{\mathbb R}^N}\int_{B(z,\sigma)}\Psi(\mu)\,\dee x
 +C\sigma^{-2}\sup_{z\in{\mathbb R}^N}\int^t_0\int_{B(z,\sigma)}u^{m-1}\Psi(u)\,\dee x\,\dee s\\
 & +C\sup_{z\in{\mathbb R}^N}\int_0^t\int_{B(z,\sigma)}\Psi(u)\,\dee x\,\dee s
 +CM_\sigma[u](t)^{p-m}\sup_{z\in{\mathbb R}^N}\int^t_0\int_{B(z,\sigma)}u^{m-1}\Psi''(u)|\nabla u|^2\,\dee x\,\dee s\
\end{split}
\end{equation}
for $t\in(0,T)$ and $\sigma>0$ if 
\begin{equation}
\label{eq:3.5}
M_\sigma[u](t):=\sup_{s\in(0,t]}\sup_{z\in\mathbb{R}^N}\sup_{r\in(0,\sigma]}
\left\{\eta(r) \Psi^{-1}\left(\,\dashint_{B(z,r)} \Psi(u(s))\,\dee y\right)\right\}\le 1.
\end{equation}
\end{lemma}
{\bf Proof.}
Let $u$ be a positive classical solution to problem~\eqref{eq:P} in ${\mathbb R}^N\times(0,T)$, where $T\in(0,\infty)$, 
and assume \eqref{eq:3.2}. 
Let $z\in{\mathbb R}^N$, $t\in(0,T)$, and $\sigma>0$. 
Let $\zeta\in C_c^\infty({\mathbb R}^N)$ be such that 
\begin{equation}
\label{eq:3.6}
\begin{split}
 & \zeta\equiv1\ \text{in}\ B(z,\sigma),\quad \zeta\equiv 0\ \text{in}\ \mathbb{R}^N\setminus B(z,2\sigma),\quad 0\le \zeta\le 1\ \text{in}\ \mathbb{R}^N,\\
 & \|\nabla\zeta\|_{L^\infty}\le 2\sigma^{-1},\quad \|\nabla^2\zeta\|_{L^\infty}\le 4\sigma^{-2}.
\end{split}
\end{equation}
Let $\Psi$ and $\eta$ be as in \eqref{eq:1.4}. Then 
\begin{equation}
\label{eq:3.7}
\begin{split}
 & \Psi'(0)>0,\quad
0<\Psi''(\xi)\preceq \xi^{-1}\Psi'(\xi)\preceq \xi^{-2}\Psi(\xi),
\quad \Psi(2\xi)\asymp \Psi(\xi),\\
& \eta'(\xi)>0,\quad \eta(2\xi)\asymp \eta(\xi),
\end{split}
\end{equation}
hold for $\xi\in(0,\infty)$.
In addition, it follows from \eqref{eq:3.7} that
\begin{equation}
\label{eq:3.8}
(\Psi^{-1})'(0)>0,
\quad (\Psi^{-1})'(\xi)>0,
\quad (\Psi^{-1})''(\xi)<0,
\quad \Psi^{-1}(2\xi)\asymp \Psi^{-1}(\xi),
\end{equation}
hold for $\xi\in(0,\infty)$.
\vspace{3pt}
\newline
\underline{Step 1.} 
Let $k\ge 1$ and $\ell\ge 1$ be large enough to be chosen later. Let $t\in(0,T)$. 
We multiply the equation 
$$
\partial_t u-\Delta u^m-u^p=0\quad\mbox{in}\quad{\mathbb R}^N\times(0,T)
$$
by $\Psi'(u\zeta^\ell)\zeta^{k+\ell}$ and integrate it in ${\mathbb R}^N\times(0,t)$. 
Since 
\begin{align*}
 & \int^t_0\int_{\mathbb{R}^N}\partial_t u\Psi'(u\zeta^\ell)\zeta^{k+\ell}\,\dee x\,\dee s\\
 &  =\int^t_0 \dfrac{\dee}{\dee s}\int_{\mathbb{R}^N}\Psi(u\zeta^\ell)\zeta^{k}\,\dee x\,\dee s
 =\int_{\mathbb{R}^N}\Psi\big(u(t)\zeta^\ell\big)\zeta^{k}\,\dee x-\int_{\mathbb{R}^N}\Psi(\mu\zeta^\ell)\zeta^{k}\,\dee x
 \end{align*}
 and
 \begin{align*}
 &-\int^t_0\int_{\mathbb{R}^N}\Delta(u^m)\Psi'(u\zeta^\ell)\zeta^{k+\ell}\,\dee x\,\dee s\\
 & =m\int^t_0\int_{\mathbb{R}^N}u^{m-1}\nabla u\cdot\bigg(\Psi''(u\zeta^\ell)\zeta^{k+\ell}(\nabla u\zeta^\ell+\ell u\zeta^{\ell-1}\nabla\zeta)
 +(k+\ell)\Psi'(u\zeta^\ell)\zeta^{k+\ell-1}\nabla \zeta\bigg)\,\dee x\,\dee s\\
 & =m\int^t_0\int_{\mathbb{R}^N}u^{m-1}\nabla u\cdot\bigg(\Psi''(u\zeta^\ell)\zeta^{k+\ell}(\nabla u\zeta^\ell+\ell u\zeta^{\ell-1}\nabla\zeta)\biggr)\,\dee x\,\dee s\\
 & \qquad
 -(k+\ell)\int^t_0\int_{\mathbb{R}^N}u^m\mbox{div}\biggr(\Psi'(u\zeta^\ell)\zeta^{k+\ell-1}\nabla \zeta\bigg)\,\dee x\,\dee s\\
 & =m\int^t_0\int_{\mathbb{R}^N}u^{m-1}\Psi''(u\zeta^\ell)|\nabla u|^2\zeta^{k+2\ell}\,\dee x\,\dee s
 +m\ell\int^t_0\int_{\mathbb{R}^N}u^m\Psi''(u\zeta^\ell)\zeta^{k+2\ell-1}\nabla u\cdot\nabla\zeta\,\dee x\,\dee s\\
 & \qquad
  -(k+\ell)\int^t_0\int_{\mathbb{R}^N} u^m\left(\Psi''(u\zeta^\ell)\zeta^{k+\ell-1}(\zeta^\ell\nabla u\cdot\nabla\zeta+\ell\zeta^{\ell-1}u|\nabla\zeta|^2)\right)\,\dee x\,\dee s\\
  & \qquad
  -\int^t_0\int_{\mathbb{R}^N} u^m\Psi'(u\zeta^\ell) \Delta\zeta^{k+\ell}\,\dee x\,\dee s\\
 & \ge\frac{m}{2}\int^t_0\int_{\mathbb{R}^N}u^{m-1}\Psi''(u\zeta^\ell)|\nabla u|^2\zeta^{k+2\ell}\,\dee x\,\dee s\\
 & \qquad
  -C\int^t_0\int_{\mathbb{R}^N} u^{m+1}\Psi''(u\zeta^\ell)\zeta^{k+2\ell-2}|\nabla\zeta|^2\,\dee x\,\dee s
   -\int^t_0\int_{\mathbb{R}^N} u^m\Psi'(u\zeta^\ell) |\Delta\zeta^{k+\ell}|\,\dee x\,\dee s,
\end{align*}
we have
\begin{equation*}
\begin{split}
	&\int_{\mathbb{R}^N}\Psi\big(u(t)\zeta^\ell\big)\zeta^{k}\,\dee x+\frac{m}{2}\int^t_0\int_{\mathbb{R}^N}u^{m-1}\Psi''(u\zeta^\ell)|\nabla u|^2\zeta^{k+2\ell}\,\dee x\,\dee s\\
	& \le \int_{\mathbb{R}^N}\Psi(\mu\zeta^\ell)\zeta^{k}\,\dee x
	+C\int^t_0\int_{\mathbb{R}^N} u^{m+1}\Psi''(u\zeta^\ell)\zeta^{k+2\ell-2}|\nabla\zeta|^2\,\dee x\,\dee s\\
	& \qquad
	+\int^t_0\int_{\mathbb{R}^N} u^m\Psi'(u\zeta^\ell) |\Delta\zeta^{k+\ell}|\,\dee x\,\dee s
	+\int^t_0\int_{\mathbb{R}^N} u^p\Psi'(u\zeta^\ell) \zeta^{k+\ell}\,\dee x\,\dee s.
\end{split}
\end{equation*}
Taking $k\ge 1$ large enough so that $k\ge (p-1)\ell$, 
by \eqref{eq:3.6} and \eqref{eq:3.7} we have
\begin{equation*}
\begin{split}
 & \int_{B(z,\sigma)}\Psi(u(t))\,\dee x+\int^t_0\int_{B(z,\sigma)}u^{m-1}\Psi''(u)|\nabla u|^2\,\dee x\,\dee s\\
 & \le C\int_{B(z,2\sigma)}\Psi(\mu)\dee x+C\sigma^{-2}\int^t_0\int_{B(z,2\sigma)}u^{m-1}\Psi(u)\,\dee x\,\dee s\\
 & \qquad\qquad
 +C\int^t_0\int_{\mathbb{R}^N}(u\zeta^\ell)^{p-1}\Psi(u\zeta^\ell)\,\dee x\,\dee s
\end{split}
\end{equation*} 
for $z\in{\mathbb R}^N$, $t\in(0,T)$, and $\sigma>0$. 
Then, by \eqref{eq:3.1} we have
\begin{equation}
\label{eq:3.9}
\begin{split}
 & \sup_{s\in(0,t]}\sup_{z\in{\mathbb R}^N}\int_{B(z,\sigma)}\Psi(u(s))\,\dee x+\sup_{z\in{\mathbb R}^N}\int^t_0\int_{B(z,\sigma)}u^{m-1}\Psi''(u)|\nabla u|^2\,\dee x\,\dee s\\
 & \le C\sup_{z\in{\mathbb R}^N}\int_{B(z,\sigma)}\Psi(\mu)\dee x
 +C\sigma^{-2}\sup_{z\in{\mathbb R}^N}\int^t_0\int_{B(z,\sigma)}u^{m-1}\Psi(u)\,\dee x\,\dee s\\
 & \qquad\qquad
 +C\sup_{z\in{\mathbb R}^N}\int^t_0\int_{\mathbb{R}^N}(u\zeta^\ell)^{p-1}\Psi(u\zeta^\ell)\,\dee x\,\dee s
\end{split}
\end{equation} 
for $t\in(0,T)$ and $\sigma>0$. 
Furthermore, it follows from $\mu\in L^\infty({\mathbb R}^N)$, \eqref{eq:3.2}, and \eqref{eq:3.9} that 
\eqref{eq:3.3} holds for $t\in(0,T)$ and $\sigma>0$. 
\vspace{3pt}
\newline
\underline{Step 2.} 
In this step we employ arguments in \cites{L} to obtain an estimate of the last term of \eqref{eq:3.9}.
Set
\begin{equation}
\label{eq:3.10}
v(x,s):=u(x,s)^{\frac{1}{\ell}}\zeta(x),\quad (x,s)\in{\mathbb R}^N\times(0,t).
\end{equation}
Then the layer cake representation (see e.g., \cite{G}*{(1.1.7)}) together with \eqref{eq:3.7} implies that
\begin{equation}
\label{eq:3.11}
\begin{split}
 & \int^t_0\int_{\mathbb{R}^N}(u\zeta^\ell)^{p-1}\Psi(u\zeta^\ell)\,\dee x\,\dee s\\
 & =\int^t_0\int_{\mathbb{R}^N} v^{\ell(p-1)}\Psi(v^\ell)\,\dee x\,\dee s\\
 & \,=\int_0^t\int_0^\infty
 \mathcal{L}^N\big(\{x\in\mathbb{R}^N\,:\,v(x, s)>\lambda\}\big)
 \frac{\dee}{\dee\lambda}\left(\lambda^{\ell(p-1)}\Psi(\lambda^\ell)\right)\,\dee\lambda\,\dee s\\
 & \,\le C\int_0^t\int_0^\infty 
 \mathcal{L}^N\big(\{x\in\mathbb{R}^N\,:\,v(x, s)>\lambda\}\big)\lambda^{\ell(p-1)-1}\Psi(\lambda^\ell)\,\dee\lambda\,\dee s\\
&  \,\le CI_\Lambda(t; z, \sigma)+C\Lambda^{\ell (p-1)}\int_0^t\int_{\mathbb{R}^N}\Psi(v^\ell)\,\dee x\,\dee s
 \end{split}
\end{equation}
for $z\in{\mathbb R}^N$, $t\in(0,T)$, $\sigma>0$, and $\Lambda\ge 0$,
where
\begin{equation}
\label{eq:3.12}
I_\Lambda(t; z, \sigma):=\int_0^t\int_\Lambda^\infty
 \mathcal{L}^N\big(\{x\in\mathbb{R}^N\,:\,v(x, s)>\lambda\}\big)\lambda^{\ell(p-1)-1}\Psi(\lambda^\ell)\,\dee\lambda\,\dee s.
\end{equation}
We obtain an estimate of $I_\Lambda(t; z, \sigma)$. 
For any $\lambda>\Lambda$ and $(x,s)\in{\mathbb R}^N\times(0,t)$, 
set
\begin{equation}
\label{eq:3.13}
v_\lambda(x, s):=
\left\{
\begin{array}{ll}
 0 & \mbox{if}\quad v(x,s)\le\displaystyle{\frac{\lambda}{2}},\vspace{3pt}\\
 v(x,s)-\displaystyle{\frac{\lambda}{2}} & \mbox{if}\quad \displaystyle{\frac{\lambda}{2}}<v(x,s)\le 2\lambda,\vspace{5pt}\\
 \displaystyle{\frac{3}{2}\lambda}  & \mbox{if}\quad v(x,s)>2\lambda.
\end{array}
\right.
\end{equation}
We claim that
\begin{equation}
\label{eq:3.14}
\begin{split}
&\dashint_{B(x,r)}v_\lambda(y, s)\,\dee y
\le\left(\frac{m_*M_\sigma[u](t)}{\eta(r)}\right)^{\frac{1}{\ell}}\\
&\qquad\qquad
 \mbox{for $x\in{\mathbb R}^N$, $r\in(0,\infty)$, $ \lambda>\Lambda$, and $s\in(0,t]$},
\end{split}
\end{equation} 
where $m_*\ge 1$ is as in \eqref{eq:3.1}.
By \eqref{eq:3.5}, \eqref{eq:3.6}, \eqref{eq:3.10}, and \eqref{eq:3.13} we apply Jensen's inequality 
to obtain 
\begin{equation*}
\begin{split}
\dashint_{B(x,r)}v_\lambda(y, s)\,\dee y
 & \le \dashint_{B(x,r)}v(y, s)\,\dee y\le \left(\,\dashint_{B(x,r)}u(y,s)\,\dee y\right)^{\frac{1}{\ell}}\\
 & \le \left(\sup_{z\in{\mathbb R}^N}\Psi^{-1}\left(\,\dashint_{B(z,r)}\Psi(u(y,s))\,\dee y\right)\right)^{\frac{1}{\ell}}
 \le\left(\frac{M_\sigma[u](t)}{\eta(r)}\right)^{\frac{1}{\ell}}
\end{split}
\end{equation*}
for $x\in{\mathbb R}^N$, $r\in(0,\sigma)$, $\lambda>\Lambda$, and $s\in(0,t]$. 
This implies that \eqref{eq:3.14} holds if $r\in(0, \sigma)$. 
On the other hand, 
since $\zeta\equiv 0$ in $\mathbb{R}^N\setminus B(z,2\sigma)$, 
by Jensen's inequality we have
$$
\dashint_{B(x,r)}v_\lambda(y, s)\,\dee y
 \le \dashint_{B(x,r)}v(y, s)\,\dee y
 \le \left(\,\frac{1}{\mathcal{L}^N(B(0,r))}\sup_{z\in{\mathbb R}^N}\,\int_{B(z,2\sigma)}u(y,s)\,\dee y\right)^{\frac{1}{\ell}}
$$
for $x\in{\mathbb R}^N$, $r\in(\sigma,\infty)$, $\lambda>\Lambda$, and $s\in(0,t]$. 
Then, by \eqref{eq:3.1} and \eqref{eq:3.5} 
we apply Jensen's inequality again to obtain 
\begin{equation*}
\begin{split}
 & \dashint_{B(x,r)}v_\lambda(y, s)\,\dee y
 \le \left(\,\frac{m_*}{\mathcal{L}^N(B(0,r))}\sup_{z\in{\mathbb R}^N}\,\int_{B(z,\sigma)}u(y,s)\,\dee y\right)^{\frac{1}{\ell}}\\
 & =\left(\,\frac{m_*\mathcal{L}^N(B(0,\sigma))}{\mathcal{L}^N(B(0,r))}\sup_{z\in{\mathbb R}^N}\,\dashint_{B(z,\sigma)}u(y,s)\,\dee y\right)^{\frac{1}{\ell}}\\
 & \le \left(\,m_*\left(\frac{\sigma}{r}\right)^N\sup_{z\in{\mathbb R}^N}\Psi^{-1}\left(\,\dashint_{B(z,\sigma)}\Psi(u(y,s))\,\dee y\right)\right)^{\frac{1}{\ell}}
  \le \left(\,m_*\left(\frac{\sigma}{r}\right)^N\frac{M_\sigma[u](t)}{\eta(\sigma)}\right)^{\frac{1}{\ell}}
\end{split}
\end{equation*} 
for $x\in{\mathbb R}^N$, $r\in(\sigma,\infty)$, $\lambda>\Lambda$, and $s\in(0,t]$. 
Since $\sigma^{-N}\eta(\sigma)\ge r^{-N}\eta(r)$ for $r\in(\sigma,\infty)$, 
we see that \eqref{eq:3.14} holds if {$r\in(\sigma, \infty)$.
Thus \eqref{eq:3.14} is valid.

Since $\eta(0)=0$ and $\lim_{\xi\to\infty}\eta(\xi)=\infty$, 
for any $\lambda>\Lambda$, we apply the intermediate value theorem
to find $r_*=r_*(\lambda)>0$ such that 
\begin{equation}
\label{eq:3.15}
\eta(r_*(\lambda))=\frac{4^\ell m_*M_\sigma[u](t)}{\lambda^\ell}.
\end{equation}
This together with \eqref{eq:3.14} implies that 
\begin{equation*}
\dashint_{B(x,r_*(\lambda))}v_\lambda(y, s)\,\dee y\le 
\left(\frac{m_*M_\sigma[u](t)}{\eta(r_*(\lambda))}\right)^{\frac{1}{\ell}}=\dfrac{\lambda}{4}
\end{equation*}
for $x\in{\mathbb R}^N$, $\lambda>\Lambda$, and $s\in(0,t]$. 
Then, by \eqref{eq:3.13} we see that
\begin{equation}
\label{eq:3.16}
\begin{split}
 & \mathcal{L}^N\big(\{x\in\mathbb{R}^N\,:\,v(x, s)>\lambda\}\big)\\
 & \le\mathcal{L}^N\left(\left\{x\in\mathbb{R}^N\,:\,v_\lambda(x, s)>\frac{\lambda}{2}\right\}\right)\\
 & \le \mathcal{L}^N\bigg(\bigg\{x\in\mathbb{R}^N\,:\,
 \bigg|\,v_\lambda(x, s)-\dashint_{B(x,r_*(\lambda))}v_\lambda(y, s)\,\dee y\,\bigg|>\dfrac{\lambda}{4}\bigg\}\bigg)\\
 & \le \dfrac{16}{\lambda^2}\int_{\mathbb{R}^N}\bigg|\,v_\lambda(x, s)-\dashint_{B(x,r_*(\lambda))}v_\lambda(y, s)\,\dee y\,\bigg|^2\,\dee x
\end{split}
\end{equation}	
for $\lambda>\Lambda$ and $s\in(0,t]$.
Furthermore, we have
\begin{equation}
\label{eq:3.17}
\begin{split}
 & \left|\,v_\lambda(x, s)-\dashint_{B(x,r_*(\lambda))}v_\lambda(y, s)\,\dee y\,\right|^2
 =\left|\,\dashint_{B(x,r_*(\lambda))}(v_\lambda(x,s)-v_\lambda(y,s))\,\dee y\,\right|^2\\
 & =\dfrac{1}{\mathcal{L}^N(B(x,r_*(\lambda)))^2}\left|\,\int_{B(x,r_*(\lambda))}
 \int_0^1\frac{\dee}{\dee\xi}v_\lambda((1-\xi)y+\xi x,s)\,\dee\xi\,\dee y\,\right|^2\\
 & \le\dfrac{1}{\mathcal{L}^N(B(x,r_*(\lambda)))^2} \left(\,\int_{B(x,r_*(\lambda))}
\int_0^1 |\nabla v_\lambda((1-\xi)y+\xi x,s)||x-y|\,\dee\xi\,\dee y\,\right)^2\\
 & \le\dfrac{r_*(\lambda)^2}{\mathcal{L}^N(B(x,r_*(\lambda)))}
 \int_{B(x,r_*(\lambda))}
 \int_0^1 |\nabla v_\lambda((1-\xi)y+\xi x,s)|^2\,\dee\xi\,\dee y\\
 & =\dfrac{r_*(\lambda)^2}{\mathcal{L}^N(B(x,r_*(\lambda)))}
 \int_0^1\int_{B(0,r_*(\lambda))}
  |\nabla v_\lambda((1-\xi)y+x,s)|^2\,\dee y\,\dee\xi.
\end{split}
\end{equation}
Then, by \eqref{eq:3.16} and \eqref{eq:3.17} we obtain
\begin{equation}
\label{eq:3.18}
\begin{split}
 & \mathcal{L}^N\big(\{x\in\mathbb{R}^N\,:\,v(x, s)>\lambda\}\big)\\
  & \le \dfrac{16r_*(\lambda)^2}{\lambda^2 \mathcal{L}^N(B(0,r_*(\lambda)))}\int^1_0\int_{B(0,r_*(\lambda))}\int_{\mathbb{R}^N}|\nabla v_\lambda((1-\xi)y+x,s)|^2\,\dee x\,\dee y\,\dee \xi\\
 & =\dfrac{16r_*(\lambda)^2}{\lambda^2}\int_{\mathbb{R}^N}|\nabla v_\lambda(x,s)|^2\,\dee x\\
 & =\dfrac{16r_*(\lambda)^2}{\lambda^2}
 \int_{\{x\in\mathbb{R}^N\,:\,\lambda/2\le v(x, s)\le 2\lambda\}}|\nabla v(x, s)|^2\,\dee x
\end{split}
\end{equation}
for $\lambda>\Lambda$ and $s\in(0,t]$.  
Since $\Psi(2\xi)\asymp \Psi(\xi)$ for $\xi\in(0,\infty)$ (see \eqref{eq:3.7}), 
we observe from \eqref{eq:3.12} and \eqref{eq:3.18} that
\begin{equation}
\label{eq:3.19}
\begin{split}
 & I_\Lambda(t;z,\sigma)\\
 & \le C\int_0^t\int_\Lambda^\infty \frac{r_*(\lambda)^2}{\lambda^2}
\left(\int_{\{x\in\mathbb{R}^N\,:\,\lambda/2\le v(s)\le 2\lambda\}}|\nabla v(s)|^2\,\dee x\right)
\lambda^{\ell(p-1)-1}\Psi(\lambda^\ell)\,\dee\lambda\,\dee s\\
& =C\int_0^t\int_\Lambda^\infty 
\lambda^{(p-1)\ell-3}r_*(\lambda)^2\Psi(\lambda^\ell)
\left(\int_{\{x\in\mathbb{R}^N\,:\,\lambda/2\le v(s)\le 2\lambda\}}|\nabla v(s)|^2\,\dee x\right)\,\dee\lambda\,\dee s\\
 & =C\int_0^t\int_{\mathbb{R}^N}|\nabla v(s)|^2\chi_{\{v(s)>\Lambda/2\}}(x)
\left(\int^{2v(s)}_{\frac{1}{2}v(s)} \lambda^{(p-1)\ell-3}r_*(\lambda)^2\Psi(\lambda^\ell)\,\dee\lambda\right)\,\dee x\,\dee s\\
 & \le C\int_0^t\int_{\mathbb{R}^N}
\left(u^{\frac{2}{\ell}-2}|\nabla u|^2\zeta^2+u^{\frac{2}{\ell}}|\nabla\zeta|^2\right)\\
& \qquad\quad
\times v^{\ell(p-1)-2}\Psi(v^\ell)\left(\sup_{\lambda\in(v(s)/2,2v(s))}r_*(\lambda)\right)^2\chi_{\{v(s)>\Lambda/2\}}(x)\,\dee x\,\dee s
\end{split}
\end{equation}
for $z\in{\mathbb R}^N$, $t\in(0,T)$, and $\sigma>0$.
\vspace{3pt}
\newline
\underline{Step 3.} 
Assume that $M_\sigma[u](t)\le 1$. 
Since $\eta$ is strictly increasing in $(0,\infty)$ (see \eqref{eq:3.7}) and 
$$
\eta^{-1}(\xi)\asymp \xi^{\frac{1}{N}}\left[\log\left(e+\frac{1}{\xi}\right)\right]^{-\frac{1}{2}},\quad \xi\in(0,\infty),
$$
it follows from \eqref{eq:3.15} that
\begin{equation}
\label{eq:3.20}
\begin{split}
\left(\sup_{\lambda\in(v(s)/2,2v(s))}r_*(\lambda)\right)^2
 &=\left[\eta^{-1}\left(\frac{8^\ell m_*M_\sigma[u](t)}{v(s)^\ell}\right)\right]^2\\
 & \le CM_\sigma[u](t)^{\frac{2}{N}}v(s)^{-\frac{2\ell}{N}}\left[\log\left(e+\frac{v(s)^\ell}{8^\ell m_*M_\sigma[u](t)}\right)\right]^{-1}\\
 & \le CM_\sigma[u](t)^{p-m}v(s)^{-\ell(p-m)}\left[\log\left(e+v(s)^\ell\right)\right]^{-1}
\end{split}
\end{equation}
for $\lambda\in(0,\infty)$. 

On the other hand, for any fixed $a>0$ and $b\in{\mathbb R}$, 
taking $L\ge e$ large enough so that 
the function $[0,\infty)\ni\xi\mapsto \xi^a[\log(L+\xi)]^b$ is increasing,
we have
\begin{equation}
\label{eq:3.21}
\xi_1^a[\log(e+\xi_1)]^b\asymp \xi_1^a[\log(L+\xi_1)]^b
\le\xi_2^a[\log(L+\xi_2)]^b\asymp \xi_2^a[\log(e+\xi_2)]^b
\end{equation}
for $\xi_1$, $\xi_2\in[0,\infty)$ with $\xi_1\le\xi_2$. 
Then, taking $\ell\ge 1$ large enough so that $m-2/\ell>0$, we observe from \eqref{eq:3.20} that
\begin{align*}
 & v(s)^{\ell(p-1)-2}\Psi(v(s)^\ell)\left(\sup_{\lambda\in(v(s)/2,2v(s))}r_*(\lambda)\right)^2\\
 & \le CM_\sigma[u](t)^{p-m}v(s)^{\ell m-2} \left[\log\left(e+v(s)^\ell\right)\right]^{\alpha-1}\\
 & \le CM_\sigma[u](t)^{p-m}u(s)^{m-\frac{2}{\ell}}\left[\log\left(e+u(s)\right)\right]^{\alpha-1}
\end{align*}
for $s\in(0,t)$. 
This together with \eqref{eq:3.1} and \eqref{eq:3.19} implies that 
\begin{equation}
\label{eq:3.22}
\begin{split}
 & I_\Lambda(t;z,\sigma)\\
 & \le CM_\sigma[u](t)^{p-m}\int_0^t\int_{B(z,2\sigma)} u^{m-3}\left[\log\left(e+u\right)\right]^{-1}\Psi(u)|\nabla u|^2\chi_{\{u(s)>\Lambda/2\}}(x)\,\dee x\,\dee s\\
 & \qquad\qquad
 +CM_\sigma[u](t)^{p-m}\sigma^{-2}\int_0^t\int_{B(z,2\sigma)}u^{m-1}\Psi(u)\,\dee x\,\dee s\\
 & \le CM_\sigma[u](t)^{p-m}\sup_{z\in{\mathbb R}^N}\int_0^t\int_{B(z,\sigma)} u^{m-3}[\log(e+u)]^{-1}\Psi(u)|\nabla u|^2\chi_{\{u(s)>\Lambda/2\}}(x)\,\dee x\,\dee s\\
 & \qquad\qquad
 +C\sigma^{-2}\sup_{z\in{\mathbb R}^N}\int_0^t\int_{B(z,\sigma)}u^{m-1}\Psi(u)\,\dee x\,\dee s
\end{split}
\end{equation}
for $z\in{\mathbb R}^N$, $t\in(0,T)$, and $\sigma>0$. 
Since 
$$
\xi^{-2}[\log(e+\xi)]^{-1}\Psi(\xi)\preceq \Psi''(\xi)\quad\mbox{for $\xi\in[1,\infty)$}, 
$$
by \eqref{eq:3.1}, \eqref{eq:3.11}, and \eqref{eq:3.22} with $\Lambda=2$ we obtain
\begin{align*}
 & \int^t_0\int_{\mathbb{R}^N}(u\zeta^\ell)^{p-1}\Psi(u\zeta^\ell)\,\dee x\,\dee s\\
 & \le C\int_0^t\int_{\mathbb{R}^N}\Psi(v^\ell)\,\dee x\,\dee s+CI_2(t; z, \sigma)\\
 & \le C\sup_{z\in{\mathbb R}^N}\int^t_0\int_{B(z,\sigma)}\Psi(u)\,\dee x\,\dee s
 +CM_\sigma[u](t)^{p-m}\sup_{z\in{\mathbb R}^N}\int_0^t\int_{B(z,\sigma)} u^{m-1}\Psi''(u)|\nabla u|^2\,\dee x\,\dee s\\
 & \qquad\qquad
 +C\sigma^{-2}\sup_{z\in{\mathbb R}^N}\int_0^t\int_{B(z,\sigma)}u^{m-1}\Psi(u)\,\dee x\,\dee s
\end{align*}
for $z\in{\mathbb R}^N$, $t\in(0,T)$, and $\sigma>0$. 
This together with \eqref{eq:3.9} implies \eqref{eq:3.4}. 
Thus Lemma~\ref{Lemma:3.1} follows.
$\Box$\vspace{3pt}

Similarly, we obtain an energy estimate of solutions to problem~\eqref{eq:P} with $p>p_m$.
\begin{lemma}
\label{Lemma:3.2}
Let $p>p_m$, $\mu\in {\mathcal L}\cap L^\infty({\mathbb R}^N)$, $T\in(0,\infty]$, and $\beta>1$. 
Let $u$ be a positive classical solution to problem~\eqref{eq:P} in ${\mathbb R}^N\times(0,T)$
satisfying \eqref{eq:3.2}.  
Then 
\begin{equation}
\label{eq:3.23}
\sup_{z\in{\mathbb R}^N}\int^t_0\int_{B(z,\sigma)}u^{m+\beta-3}|\nabla u|^2\,\dee x\,\dee s<\infty
\end{equation}
for $t\in(0,T)$ and $\sigma>0$. 
Furthermore, there exists $C=C(N,m,p,\beta)>0$ such that 
\begin{equation}
\label{eq:3.24}
\begin{split} 
 & \sup_{s\in(0,t]}\sup_{z\in{\mathbb R}^N}\int_{B(z,\sigma)}u(s)^\beta\,\dee x+\sup_{z\in{\mathbb R}^N}\int^t_0\int_{B(z,\sigma)}u^{m+\beta-3}|\nabla u|^2\,\dee x\,\dee s\\
 & \le C\sup_{z\in{\mathbb R}^N}\int_{B(z,\sigma)}\mu^\beta\,\dee x\\
 & \qquad
  +C\left(1+\sup_{s\in(0,t)}|||u(s)|||^{p-m}_{\frac{N(p-m)}{2},\beta;\sigma}\right)
 \sigma^{-2}\sup_{z\in{\mathbb R}^N}\int^t_0\int_{B(z,\sigma)}u^{m+\beta-1}\,\dee x\,\dee s\\
 & \qquad
 +C\left(\sup_{s\in(0,t)}|||u(s)|||_{\frac{N(p-m)}{2},\beta;\sigma}\right)^{p-m} 
 \sup_{z\in{\mathbb R}^N}\int^t_0\int_{B(z,\sigma)}u^{m+\beta-3}|\nabla u|^2\,\dee x\,\dee s
\end{split}
\end{equation}
for $t\in(0,T)$ and $\sigma>0$. 
\end{lemma}
{\bf Proof.} 
Setting 
$\Psi(\xi)=\xi^\beta$ and $\eta(\xi)=\xi^{\frac{2}{p-m}}$ for $\xi\in[0,\infty)$, 
we apply the same arguments as in the proof of Lemma~\ref{Lemma:3.1}. 
Then we have
\begin{equation}
\label{eq:3.25}
\begin{split}
 & \sup_{z\in{\mathbb R}^N}\int_{B(z,\sigma)}u(t)^\beta\,\dee x+\sup_{z\in{\mathbb R}^N}\int^t_0\int_{B(z,\sigma)}u^{m+\beta-3}|\nabla u|^2\,\dee x\,\dee s\\
 & \le C\sup_{z\in{\mathbb R}^N}\int_{B(z,\sigma)}\mu^\beta\dee x
 +C\sigma^{-2}\sup_{z\in{\mathbb R}^N}\int^t_0\int_{B(z,\sigma)}u^{m+\beta-1}\,\dee x\,\dee s\\
 & \qquad\qquad
 +C\sup_{z\in{\mathbb R}^N}\int^t_0\int_{\mathbb{R}^N}(u\zeta^\ell)^{p+\beta-1}\,\dee x\,\dee s
\end{split}
\end{equation} 
for $t\in(0,T)$ and $\sigma>0$, instead of \eqref{eq:3.9}. 
This together with \eqref{eq:3.2} implies \eqref{eq:3.23}. 

Let $I_\Lambda$, $v$, and $r_*$ be as in the proof of Lemma~\ref{Lemma:3.1}. 
It follows from \eqref{eq:3.15} that
$$
r_*(\lambda)^{\frac{2}{p-m}}=\frac{4^\ell m_*M_\sigma[u](t)}{\lambda^\ell}
$$
and hence
$$
\left(\sup_{\lambda\in(v(s)/2,2v(s))}r_*(\lambda)\right)^2\le CM_\sigma[u](t)^{p-m}(u\zeta^\ell)^{-(p-m)}.
$$
Then, taking $\ell\ge 1$ large enough if necessary, 
by \eqref{eq:3.19} with $\Lambda=0$ we obtain 
\begin{equation}
\label{eq:3.26}
\begin{split}
 & \int^t_0\int_{\mathbb{R}^N}(u\zeta^\ell)^{p-1}\Psi(u\zeta^\ell)\,\dee x\,\dee s =I_0(t;z,\sigma)\\
 & \le CM_\sigma[u](t)^{p-m}\int_0^t\int_{B(z,2\sigma)}u^{m-3}\Psi(u)|\nabla u|^2\,\dee x\,\dee s\\
 & \qquad
 +CM_\sigma[u](t)^{p-m}\sigma^{-2}\int_0^t\int_{B(z,2\sigma)}u^{m-1}\Psi(u)\,\dee x\,\dee s\\
 & \le CM_\sigma[u](t)^{p-m}\sup_{z\in{\mathbb R}^N}\int_0^t\int_{B(z,\sigma)}u^{m-3}\Psi(u)|\nabla u|^2\,\dee x\,\dee s\\
  & \qquad
 +CM_\sigma[u](t)^{p-m}\sigma^{-2}\sup_{z\in{\mathbb R}^N}\int_0^t\int_{B(z,\sigma)}u^{m-1}\Psi(u)\,\dee x\,\dee s
\end{split}
\end{equation}
for $z\in{\mathbb R}^N$, $t\in(0,T)$, and $\sigma>0$.
Since
$$
M_\sigma[u](t)=\sup_{s\in(0,t)}|||u(s)|||_{\frac{N(p-m)}{2},\beta;\sigma}
$$
for $t\in(0,T)$, 
by \eqref{eq:3.25} and \eqref{eq:3.26} we obtain \eqref{eq:3.24}. 
Thus Lemma~\ref{Lemma:3.2} follows.
$\Box$
\vspace{5pt}

At the end of this section we recall 
decay estimates of solutions to problem~\eqref{eq:P}. 
See \cite{ADi}*{Proposition~7.1}. (See also \cite{Sato}*{Proposition~3.2}.)
\begin{lemma}
\label{Lemma:3.3}
Let $u$ be a solution to problem~\eqref{eq:P} in ${\mathbb R}^N\times(0,T)$ for some $T\in(0,\infty)$. 
Let $r\ge 1$. Then there exists $C=C(N,m,p,r)>0$ such that 
\begin{equation*}
\|u(t)\|_{L^\infty({\mathbb R}^N)}\le Ct^{-\frac{N+2}{\kappa_r}}\sup_{z\in{\mathbb R}^N}\left(\int_{t/2}^t\int_{B(z,2\sigma)} u^r\,\dee x\,\dee s\right)^{\frac{2}{\kappa_r}}
\end{equation*}
for $t\in(0,T_*)$,
where $\kappa_r:=N(m-1)+2r$ and
$$
T_*:=\sup\left\{\tau\in(0,T)\,:\,\sigma^{-2}\|u(s)\|^{m-1}_{L^\infty({\mathbb R}^N)}+\|u(s)\|_{L^\infty({\mathbb R}^N)}^{p-1}\le 2s^{-1}\mbox{ for $s\in(0,\tau)$}\right\}.
$$
\end{lemma}
\section{Proof of Theorem~\ref{Theorem:1.2}}
In this section we study sufficient conditions for the existence of solutions to problem~\eqref{eq:P} with $p=p_m$, 
and prove Theorem~\ref{Theorem:1.2}.

Let $\Psi$ and $\eta$ be as in Theorem~\ref{Theorem:1.2}. 
Define a $C^1$-function $\gamma$ in $[0,1]$ by 
\begin{equation}
\label{eq:4.1}
\int^{\gamma(\xi)}_0 s\eta(s)^{m-1}\,\dee s=C_\eta \xi\quad\mbox{for $\xi\in[0,1]$},
\quad\mbox{where $\displaystyle{C_\eta:=\int^1_0 s\eta(s)^{m-1}\,\dee s}$}.
\end{equation}
Then $\gamma'>0$ in $[0,1]$, $\gamma(0)=0$, and $\gamma(1)=1$. 
Since 
\begin{align*}
C_\eta\xi=\int^{\gamma(\xi)}_0 s\eta(s)^{m-1}\,\dee s & =\int^{\gamma(\xi)}_0 s^{1+N(m-1)}\left[\log\left(e+\frac{1}{s}\right)\right]^{\frac{N(m-1)}{2}}\,\dee s\\
 & \asymp \gamma(\xi)^{2+N(m-1)}\left[\log\left(e+\frac{1}{\gamma(\xi)}\right)\right]^{\frac{N(m-1)}{2}}
=\gamma(\xi)^2\eta(\gamma(\xi))^{m-1}
\end{align*}
for $\xi\in[0,1]$, we see that 
\begin{equation}
\label{eq:4.2}
\gamma(\xi)^2\eta(\gamma(\xi))^{m-1}\asymp\xi\quad\mbox{for}\quad \xi\in[0,1]. 
\end{equation}
We prove a lemma on the function $\gamma$. 
\begin{lemma}
\label{Lemma:4.1}
Let $p=p_m$. 
The function $\gamma$ defined by \eqref{eq:4.1} satisfies the following properties:
\begin{align}
\label{eq:4.3}
 & \gamma(\xi)\asymp\gamma\left(\frac{\xi}{2}\right),
 \quad \eta(\gamma(\xi))\asymp\eta\left(\gamma\left(\frac{\xi}{2}\right)\right),\\
\label{eq:4.4}
 & \eta(\gamma(\xi))\succeq \xi^{\frac{1}{p-1}},\\
\label{eq:4.5} 
 & \int_0^\xi \eta(\gamma(s))^{-(m-1)}\,\dee s\asymp \gamma(\xi)^2,
\end{align}
for $\xi\in[0,1]$. 
\end{lemma}
{\bf Proof.}
We prove \eqref{eq:4.3}. 
Taking $k\ge 1$ large enough, by \eqref{eq:4.1} and the monotonicity of $\eta$ 
we obtain 
\begin{align*}
\int^{k^{-1}\gamma(\xi)}_0 s\eta(s)^{m-1}\,\dee s
 & =k^{-2}\int^{\gamma(\xi)}_0 s\eta(k^{-1}s)^{m-1}\,\dee s\\
 & \le Ck^{-2}\int^{\gamma(\xi)}_0 s\eta(s)^{m-1}\,\dee s
=Ck^{-2} C_\eta \xi\le\frac{\xi}{2}C_\eta.
\end{align*}
Then we observe from \eqref{eq:4.1} that 
$$
k^{-1}\gamma(\xi)\le\gamma\left(\frac{\xi}{2}\right),\quad \xi\in[0,1].
$$
This together with the monotonicity of $\gamma$ and the relation that $\eta(2\xi)\asymp\eta(\xi)$ for $\xi\in(0,\infty)$ (see \eqref{eq:3.7}) 
implies \eqref{eq:4.3}. 

We prove \eqref{eq:4.4}. 
Since $\eta(\xi)\ge \xi^N=\xi^{2/(p-m)}$ for $\xi\in[0,\infty)$, 
it follows from \eqref{eq:4.1} that
$$
C_\eta \xi\ge C\int_0^{\gamma(\xi)} s^{1+\frac{2(m-1)}{p-m}}\,\dee s= C\gamma(\xi)^{2+\frac{2(m-1)}{p-m}}=C\gamma(\xi)^{\frac{2(p-1)}{p-m}}
$$
for $\xi\in[0,1]$. 
This together with \eqref{eq:4.2} implies that
$$
\eta(\gamma(\xi))\asymp \left(\xi\gamma(\xi)^{-2}\right)^{\frac{1}{m-1}}\succeq \xi^{\frac{1}{m-1}} \left(\xi^{-\frac{p-m}{2(p-1)}}\right)^{\frac{2}{m-1}}
=\xi^{\frac{1}{p-1}}\quad\mbox{for}\quad \xi\in[0,1].
$$
Thus \eqref{eq:4.4} holds. 

It remains to prove \eqref{eq:4.5}. 
It follows from \eqref{eq:4.1} that
$$
\gamma(\xi)\gamma(\xi)'\eta(\gamma(\xi))^{m-1}=C_\eta\quad\mbox{for $\xi\in[0,1]$}, 
$$
that is, 
$$
\frac{\dee}{\dee\xi}\gamma(\xi)^2=2C_\eta \eta(\gamma(\xi))^{-(m-1)}.
$$
This together with $\gamma(0)=0$ implies that
$$
\gamma(\xi)^2=2C_\eta\int_0^\xi \eta(\gamma(s))^{-(m-1)}\,\dee s,\quad\xi\in[0,1].
$$
Thus \eqref{eq:4.5} holds, and the proof of Lemma~\ref{Lemma:4.1} is complete.
$\Box$
\vspace{5pt}

Now we are ready to prove Theorem~\ref{Theorem:1.2}.
\vspace{3pt}
\newline
{\bf Proof of Theorem~\ref{Theorem:1.2}.}
By Remark~\ref{Remark:1.1}-(1) 
it suffices to consider the case of $T=1$. 
Let $p=p_m$ and $\mu\in{\mathcal L}$. 
Let $\Psi$ and $\eta$ be as in Theorem~\ref{Theorem:1.2}.
Let $\epsilon_5>0$ be small enough, and assume \eqref{eq:1.5}.
For any $i$, $j=1,2,\dots$, 
let $u_{ij}$ be a solution to problem~\eqref{eq:P} with initial data $\mu$ replaced by
\begin{equation}
\label{eq:4.6}
\mu_{ij}(x):=\min\{\mu(x),i\}+j^{-1},\quad \mbox{a.a.~$x\in{\mathbb R}^N$}. 
\end{equation}
By arguments in \cite{ADi} we find a unique classical solution $u_{ij}$ to problem~\eqref{eq:P} in ${\mathbb R}^N\times(0,T_{ij})$ 
such that 
\begin{align}
\label{eq:4.7}
 & u_{ij}(x,t)\ge j^{-1}\quad\mbox{for $(x,t)\in{\mathbb R}^N\times(0,T_{ij})$},\\
\label{eq:4.8}
 & \sup_{t\in(0,T)}\|u_{ij}(t)\|_{L^\infty({\mathbb R}^N)}<\infty\quad\mbox{for $T\in(0,T_{ij})$},\\
\label{eq:4.9}
 & \limsup_{t\nearrow T_{ij}}\|u_{ij}(t)\|_{L^\infty({\mathbb R}^N)}=\infty,
\end{align}
where $T_{ij}$ is the maximal existence time of $u_{ij}$. 
Since $\Psi^{-1}$ is Lipschitz continuous on $[0, \infty)$ (see \eqref{eq:3.8}), 
by \eqref{eq:3.7} and \eqref{eq:3.8} we find $j_*=j_*(\epsilon_5)\in\mathbb{N}$ such that
\begin{equation}
\label{eq:4.10}
\begin{split}
|||\mu_{ij}|||_{\eta,\Psi;1}
 & \le\sup_{z\in{\mathbb R}^N}\sup_{\sigma\in(0,1]}
\eta(\sigma)\Psi^{-1}\left(\,\dashint_{B(z,\sigma)}\Psi(\mu+j^{-1})\,\dee x\right)\\
 & \le C\sup_{z\in{\mathbb R}^N}\sup_{\sigma\in(0,1]}
\eta(\sigma)\Psi^{-1}\left(\,\dashint_{B(z,\sigma)}\Psi(\mu)\,\dee x+j^{-1}\right)\\
 & \le C\sup_{z\in{\mathbb R}^N}\sup_{\sigma\in(0,1]}
\eta(\sigma)\Psi^{-1}\left(\,\dashint_{B(z,\sigma)}\Psi(\mu)\,\dee x\right)+Cj^{-1}\le C\epsilon_5
\end{split}
\end{equation}
for $i\ge 1$ and $j\ge j_*$. 
\vspace{3pt}
\newline 
\underline{Step 1.}
Let $\delta_1$, $\delta_2\in(0,1)$. Set
\begin{align}
\label{eq:4.11}
T^1_{ij}:= & \,\sup\left\{T\in(0,T_{ij})\,:\, \sup_{t\in(0,T)}\,
\eta(\gamma(t))\|u_{ij}(t)\|_{L^\infty({\mathbb R}^N)}\le \delta_1\right\},\\
\label{eq:4.12}
T^2_{ij}:= & \,\sup\left\{T\in(0,T_{ij})\,:\,\sup_{t\in(0,T)}|||u_{ij}(t)|||_{\eta,\Psi;1}\le\delta_2\right\}.
\end{align}
Under suitable choices of $\delta_1$ and $\delta_2$, 
taking $\epsilon_5>0$ small enough if necessary, 
we show that 
$$
T^*_{ij}:=\min\{T^1_{ij},T^2_{ij},1\}=1\quad\mbox{for $i\ge 1$ and $j\ge j_*$}. 
$$ 
In the proof of Theorem~\ref{Theorem:1.2}, the constants $C$ are independent of $i\ge 1$ and $j\ge j_*$.

We first show that $T^*_{ij}>0$.
Since $T^1_{ij}>0$ immediately follows from \eqref{eq:4.8}, it suffices to show that $T^2_{ij}>0$.
Set $c_{ij}:=\sup_{s\in(0, T_{ij}/2)}\|u_{ij}(s)\|_{L^\infty({\mathbb R}^N)}<\infty$ and take $t_{ij}\in(0, \min\{T_{ij}/2, 1\})$ satisfying
\begin{equation}
\label{eq:t_{ij}}
\eta(t^{1/4}_{ij})c_{ij}\le \dfrac{\delta_2}{2}\quad\text{and} \quad t^{1/2}_{ij}(c_{ij}^{m-1}+c_{ij}^{p-1})\le 1
\end{equation}
for $i\ge 1$ and $j\ge j_*$.
Then we have
\begin{equation}
\label{eq:T^2}
\begin{split}
&\sup_{s\in(0,t_{ij})}\sup_{z\in{\mathbb R}^N}\sup_{\sigma\in(0, t^{1/4}_{ij}]}
\eta(\sigma)\Psi^{-1}\left(\,\dashint_{B(z,\sigma)} \Psi(u_{ij}(s))\,\dee y\right)\\
&\le \sup_{s\in(0,t_{ij})}\sup_{z\in{\mathbb R}^N}\sup_{\sigma\in(0, t^{1/4}_{ij}]}
\eta(\sigma)\Psi^{-1}\left(\Psi\left(\|u_{ij}(s)\|_{L^\infty({\mathbb R}^N)}\right)\right)\le \eta(t^{1/4}_{ij})c_{ij}\le \dfrac{\delta_2}{2}
\end{split}
\end{equation}
for $i\ge 1$ and $j\ge j_*$.
On the other hand, similarly to the argument in Step.1 in the proof of Lemma~\ref{Lemma:3.1} (see \eqref{eq:3.9}), 
it follows from \eqref{eq:t_{ij}} that
\begin{equation*}
\begin{split}
 & \sup_{s\in(0,t]}\sup_{z\in{\mathbb R}^N}\int_{B(z,\sigma)}\Psi(u_{ij}(s))\,\dee x\\
 &\le C\sup_{z\in{\mathbb R}^N}\int_{B(z,\sigma)}\Psi(\mu_{ij})\dee x+C(\sigma^{-2}c_{ij}^{m-1}+c_{ij}^{p-1})\sup_{z\in{\mathbb R}^N}\int^t_0\int_{B(z,\sigma)}\Psi(u_{ij})\,\dee x\,\dee s\\
 &\le C\sup_{z\in{\mathbb R}^N}\int_{B(z,\sigma)}\Psi(\mu_{ij})\dee x+Ct_{ij}^{-1}\sup_{z\in{\mathbb R}^N}\int^t_0\int_{B(z,\sigma)}\Psi(u_{ij})\,\dee x\,\dee s
\end{split}
\end{equation*}
for $t\in(0, t_{ij})$, $\sigma\in(t^{1/4}_{ij},1]$, $i\ge 1$, and $j\ge j_*$.
This together with Gronwall's inequality implies that
$$
 \sup_{s\in(0,t]}\sup_{z\in{\mathbb R}^N}\dashint_{B(z,\sigma)}\Psi(u_{ij}(s))\,\dee x\le Ce^{Ct^{-1}_{ij}t}\sup_{z\in{\mathbb R}^N}\dashint_{B(z,\sigma)}\Psi(\mu_{ij})\dee x\le C\sup_{z\in{\mathbb R}^N}\dashint_{B(z,\sigma)}\Psi(\mu_{ij})\dee x
$$
for $t\in(0, t_{ij})$, $\sigma\in(t^{1/4}_{ij},1]$, $i\ge 1$, and $j\ge j_*$.
Since $\Psi^{-1}(2\xi)\asymp \Psi^{-1}(\xi)$ for $\xi\in(0,\infty)$ (see \eqref{eq:3.7}), it follows from \eqref{eq:4.10} that
$$
\sup_{s\in(0,t_{ij}]}\sup_{z\in{\mathbb R}^N}\sup_{\sigma\in(t^{1/4}_{ij}, 1]}\eta(\sigma)\Psi^{-1}\left(\dashint_{B(z,\sigma)}\Psi(u_{ij}(s))\,\dee x\right)\le C\epsilon_{5}
$$
for $i\ge 1$ and $j\ge j_*$.
This together with \eqref{eq:T^2} implies that, taking $\epsilon_5>0$ small enough if necessary, we obtain $T^2_{ij}\ge t_{ij}>0$.

By Lemma~\ref{Lemma:3.1} and $\delta_2\in(0,1)$ (see \eqref{eq:4.12}) 
we obtain
\begin{equation*}
\begin{split} 
 & \sup_{s\in(0,t]}\sup_{z\in{\mathbb R}^N}\int_{B(z,\sigma)}\Psi(u_{ij}(s))\,\dee x+\sup_{z\in{\mathbb R}^N}\int^t_0\int_{B(z,\sigma)}u_{ij}^{m-1}\Psi''(u_{ij})|\nabla u_{ij}|^2\,\dee x\,\dee s\\
 & \le C\sup_{z\in{\mathbb R}^N}\int_{B(z,\sigma)}\Psi(\mu_{ij})\,\dee x
 +C\sigma^{-2}\sup_{z\in{\mathbb R}^N}\int^t_0\int_{B(z,\sigma)}u_{ij}^{m-1}\Psi(u_{ij})\,\dee x\,\dee s\\
 & +C\sup_{z\in{\mathbb R}^N}\int_0^t\int_{B(z,\sigma)}\Psi(u_{ij})\,\dee x\,\dee s
 +C\delta_2^{p-m}\sup_{z\in{\mathbb R}^N}\int^t_0\int_{B(z,\sigma)}u_{ij}^{m-1}\Psi''(u_{ij})|\nabla u_{ij}|^2\,\dee x\,\dee s,
\end{split}
\end{equation*}
for $t\in(0,T_{ij}^*)$, $\sigma\in(0,1]$, $i\ge 1$, and $j\ge j_*$.
Then, taking $\delta_2>0$ small enough if necessary, we have
\begin{equation*}
\begin{split} 
 & \sup_{z\in{\mathbb R}^N}\int_{B(z,\sigma)}\Psi(u_{ij}(t))\,\dee x\\
 & \le C\sup_{z\in{\mathbb R}^N}\int_{B(z,\sigma)}\Psi(\mu_{ij})\,\dee x
 +C\sigma^{-2}\sup_{z\in{\mathbb R}^N}\int^t_0\int_{B(z,\sigma)}u_{ij}^{m-1}\Psi(u_{ij})\,\dee x\,\dee s\\
 & +C\sup_{z\in{\mathbb R}^N}\int_0^t\int_{B(z,\sigma)}\Psi(u_{ij})\,\dee x\,\dee s
\end{split}
\end{equation*}
for $t\in(0,T_{ij}^*)$, $\sigma\in(0,1]$, $i\ge 1$, and $j\ge j_*$.
This together with \eqref{eq:4.11} implies that 
\begin{equation*}
\begin{split} 
X_{ij}(t)
  & \le C\sup_{z\in{\mathbb R}^N}\,\dashint_{B(z,\sigma)}\Psi(\mu_{ij})\,\dee x
+C\sigma^{-2}\int_0^t \|u_{ij}(s)\|_{L^\infty({\mathbb R}^N)}^{m-1}X_{ij}(s)\,\dee s+C\int_0^t X_{ij}(s)\,\dee s\\
& \le C\sup_{z\in{\mathbb R}^N}\,\dashint_{B(z,\sigma)}\Psi(\mu_{ij})\,\dee x
+C\sigma^{-2}\int_0^t \eta(\gamma(s))^{-(m-1)}X_{ij}(s)\,\dee s
+C\int_0^t X_{ij}(s)\,\dee s
\end{split}
\end{equation*}
for $t\in(0,T_{ij}^*)$, $\sigma\in(0,1]$, $i\ge 1$, and $j\ge j_*$, where
$$
X_{ij}(t):=\sup_{z\in{\mathbb R}^N}\dashint_{B(z,\sigma)}\Psi(u_{ij}(t))\,\dee x.
$$
Then Gronwall's inequality together with \eqref{eq:4.5} implies that
\begin{align*}
X_{ij}(t) & \le C\exp\left(C\sigma^{-2}\int_0^t \eta(\gamma(s))^{-(m-1)}\,\dee s+Ct\right)
\sup_{z\in{\mathbb R}^N}\dashint_{B(z,\sigma)}\Psi(\mu_{ij})\dee x\\
 & \le C\exp\left(C\sigma^{-2}\gamma(t)^2\right)
\sup_{z\in{\mathbb R}^N}\dashint_{B(z,\sigma)}\Psi(\mu_{ij})\dee x
\le C\sup_{z\in{\mathbb R}^N}\dashint_{B(z,\sigma)}\Psi(\mu_{ij})\dee x
\end{align*}
for $t\in(0,T^*_{ij})$, $\sigma\in[\gamma(t),1]$, $i\ge 1$, and $j\ge j_*$.
Then, since $\Psi^{-1}(2\xi)\asymp \Psi^{-1}(\xi)$ for $\xi\in(0,\infty)$ (see \eqref{eq:3.7}), 
we obtain
$$
\sup_{z\in{\mathbb R}^N}\Psi^{-1}\left(\dashint_{B(z,\sigma)}\Psi(u_{ij}(t))\,\dee x\right)
\le C\sup_{z\in{\mathbb R}^N}\Psi^{-1}\left(\dashint_{B(z,\sigma)}\Psi(\mu_{ij})\,\dee x\right)
$$
for $t\in(0,T_{ij}^*)$, $\sigma\in[\gamma(t),1]$, $i\ge 1$, and $j\ge j_*$.
Therefore, thanks to \eqref{eq:4.10}, we obtain 
\begin{equation}
\label{eq:4.13}
\begin{split}
 & \sup_{s\in(0, t)}\sup_{z\in{\mathbb R}^N}\sup_{\sigma\in[\gamma(s),1]}\eta(\sigma)\Psi^{-1}\left(\,\dashint_{B(z,\sigma)}\Psi(u_{ij}(s))\,\dee x\right)\\
 & \le C\sup_{z\in{\mathbb R}^N}\sup_{\sigma\in (0,1]}\eta(\sigma)\Psi^{-1}\left(\,\dashint_{B(z,\sigma)}\Psi(\mu_{ij})\,\dee x\right)
  = C|||\mu_{ij}|||_{\eta,\Psi;1}\le C\epsilon_5
\end{split}
\end{equation}
for $t\in(0,T_{ij}^*)$, $i\ge 1$, and $j\ge j_*$.

On the other hand, taking $\delta_1>0$ small enough, 
by \eqref{eq:4.2}, \eqref{eq:4.4}, and \eqref{eq:4.11} we have 
\begin{align*}
 & \left(\frac{\gamma(s)}{2}\right)^{-2}\|u_{ij}(s)\|^{m-1}_{L^\infty({\mathbb R}^N)}+\|u_{ij}(s)\|_{L^\infty({\mathbb R}^N)}^{p-1}\\
 & \le 4\delta_1^{m-1} \gamma(s)^{-2}\eta(\gamma(s))^{-(m-1)}+\delta_1^{p-1}\eta(\gamma(s))^{-(p-1)}\\
 & \le C\delta_1^{m-1}s^{-1}+C\delta_1^{p-1}s^{-1}\le 2s^{-1}
\end{align*}
for $s\in(0,T_{ij}^*)$, $i\ge 1$, and $j\ge j_*$.
Then, by Lemma~\ref{Lemma:3.3} with $\sigma=\gamma(s)/2$ and $r=1$ we have
\begin{equation*}
\begin{split}
\|u_{ij}(t)\|_{L^\infty({\mathbb R}^N)}
 & \le Ct^{-\frac{N+2}{\kappa_1}}\left(\sup_{z\in{\mathbb R}^N}\int_{t/2}^t\int_{B(z,\gamma(t))} u_{ij}\,dx\,ds\right)^{\frac{2}{\kappa_1}}\\
 & \le Ct^{-\frac{N}{\kappa_1}}\gamma(t)^{\frac{2N}{\kappa_1}}\sup_{s\in(t/2,t)}
 \sup_{z\in{\mathbb R}^N}\left(\,\dashint_{B(z,\gamma(t))} u_{ij}(s)\,\dee x\right)^{\frac{2}{\kappa_1}}
 \end{split}
\end{equation*}
for $t\in(0,T_{ij}^*)$, $i\ge 1$, and $j\ge j_*$.
This together with Jensen's inequality, \eqref{eq:4.2}, and \eqref{eq:4.3} implies that
\begin{equation}
\label{eq:4.14}
\begin{split}
  & \eta(\gamma(t))\|u_{ij}(t)\|_{L^\infty({\mathbb R}^N)}\\
 & \le Ct^{-\frac{N}{\kappa_1}}\gamma(t)^{\frac{2N}{\kappa_1}}\eta(\gamma(t))^{1-\frac{2}{\kappa_1}}\sup_{s\in(t/2,t)}
 \sup_{z\in{\mathbb R}^N}\left[\eta(\gamma(t))\Psi^{-1}\left(\,\dashint_{B(z,\gamma(t))} \Psi(u_{ij}(s))\,\dee x\right)\right]^{\frac{2}{\kappa_1}}\\
  & \le C\left(t^{-1}\gamma(t)^2\eta(\gamma(t))^{m-1}\right)^\frac{N}{\kappa_1}\sup_{s\in(0,t)}\sup_{z\in{\mathbb R}^N}\sup_{\sigma\in[\gamma(s),1]}
  \left[\eta(\sigma)\Psi^{-1}\left(\,\dashint_{B(z,\sigma)} \Psi(u_{ij}(s))\,\dee x\right)\right]^{\frac{2}{\kappa_1}}\\ 
  & \le C\sup_{s\in(0,t)}\sup_{z\in{\mathbb R}^N}\sup_{\sigma\in[\gamma(s),1]}
  \left[\eta(\sigma)\Psi^{-1}\left(\,\dashint_{B(z,\sigma)} \Psi(u_{ij}(s))\,\dee x\right)\right]^{\frac{2}{\kappa_1}}
 \end{split}
\end{equation}
for $t\in(0,T_{ij}^*)$, $i\ge 1$, and $j\ge j_*$.
Therefore, 
by \eqref{eq:4.13} and \eqref{eq:4.14}, 
taking $\epsilon_5>0$ small enough if necessary, 
we obtain 
\begin{equation}
\label{eq:4.15}
\eta(\gamma(t))\|u_{ij}(t)\|_{L^\infty({\mathbb R}^N)}\le C|||\mu_{ij}|||_{\eta,\Psi;1}^{\frac{2}{\kappa_1}}\le C\epsilon_5^{\frac{2}{\kappa_1}}\le\frac{\delta_1}{2}
\end{equation}
for $t\in(0,T^*_{ij})$, $i\ge 1$, and $j\ge j_*$.
Furthermore, we observe from \eqref{eq:4.15} that 
\begin{equation}
\label{eq:4.16}
\begin{split}
 & \sup_{s\in(0,t)}\sup_{z\in{\mathbb R}^N}\sup_{\sigma\in(0,\gamma(s))}
\eta(\sigma)\Psi^{-1}\left(\,\dashint_{B(z,\sigma)} \Psi(u_{ij}(s))\,\dee y\right)\\
 & \le\sup_{s\in(0,t)}\sup_{\sigma\in(0,\gamma(s))}
\eta(\sigma)\Psi^{-1}\left(\Psi\left(\|u_{ij}(s)\|_{L^\infty({\mathbb R}^N)}\right)\right)
=\sup_{s\in(0,t)}\eta(\gamma(s))\|u_{ij}(s)\|_{L^\infty({\mathbb R}^N)}\\
 & \le C|||\mu_{ij}|||_{\eta,\Psi;1}^{\frac{2}{\kappa_1}}
  \le C\epsilon_5^{\frac{2}{\kappa_1}}
 \end{split}
\end{equation}
for $t\in(0,T_{ij}^*)$, $i\ge 1$, and $j\ge j_*$.
Combining \eqref{eq:4.13} and \eqref{eq:4.16} 
and taking $\epsilon_5>0$ small enough if necessary, 
we obtain 
\begin{equation}
\label{eq:4.17}
\sup_{s\in(0,t)}|||u_{ij}(s)|||_{\eta,\Psi;1}\le C\epsilon_5+C\epsilon_5^{\frac{2}{\kappa_1}}\le\frac{\delta_2}{2}
\end{equation}
for $t\in(0,T_{ij}^*)$, $i\ge 1$, and $j\ge j_*$.
Then, thanks to \eqref{eq:4.15} and \eqref{eq:4.17}, 
by the definition of $T_{ij}^*$ we see that $T_{ij}^*=1$ for $i\ge 1$ and $j\ge j_*$. 
Furthermore, 
repeating the arguments in \eqref{eq:4.10}, \eqref{eq:4.13}, \eqref{eq:4.15}, and \eqref{eq:4.16}, we obtain 
\begin{equation}
\label{eq:4.18}
\sup_{s\in(0,1)}\eta(\gamma(s))\|u_{ij}(s)\|_{L^\infty({\mathbb R}^N)}+
\sup_{s\in(0,1)}|||u_{ij}(s)|||_{\eta,\Psi;1}\le C\left(|||\mu|||_{\eta,\Psi;1}+j^{-1}\right)^{\frac{2}{\kappa_1}}
\end{equation}
for $i\ge 1$ and $j\ge j_*$. 
On the other hand, it follows from \eqref{eq:4.2} that
\begin{align}
\label{eq:4.19}
 & \eta(\gamma(\xi))\asymp \xi^{\frac{1}{m-1}}\gamma(\xi)^{-\frac{2}{m-1}},\\
\label{eq:4.20}
 & \gamma(\xi)^{N(m-1)+2}\left[\log\left(e+\frac{1}{\gamma(\xi)}\right)\right]^{\frac{N(m-1)}{2}}\asymp\xi,
\end{align}
for $\xi\in(0,1)$. Then, by \eqref{eq:4.20} we have
$$
\gamma(\xi)\asymp \xi^{\frac{1}{N(m-1)+2}}\left[\log\left(e+\frac{1}{\xi}\right)\right]^{-\frac{m-1}{2}\frac{N}{N(m-1)+2}}
\quad\mbox{for $\xi\in(0,1)$}.
$$
This together with \eqref{eq:4.19} implies that 
\begin{equation}
\label{eq:4.21}
\begin{split}
\eta(\gamma(\xi))& \asymp \xi^{\frac{1}{m-1}\left(1-\frac{2}{N(m-1)+2}\right)}\left[\log\left(e+\frac{1}{\xi}\right)\right]^{\frac{N}{N(m-1)+2}}\\
 & = \xi^{\frac{N}{N(m-1)+2}}\left[\log\left(e+\frac{1}{\xi}\right)\right]^{\frac{N}{N(m-1)+2}}
 = \xi^{\frac{1}{p-1}}\left[\log\left(e+\frac{1}{\xi}\right)\right]^{\frac{1}{p-1}}
\end{split}
\end{equation}
for $\xi\in(0,1)$. 
We deduce from \eqref{eq:4.18} and \eqref{eq:4.21} that
\begin{equation}
\label{eq:4.22}
\begin{split}
&\sup_{s\in(0,1)}s^{\frac{1}{p-1}}\left[\log\left(e+\frac{1}{s}\right)\right]^{\frac{1}{p-1}}\|u_{ij}(s)\|_{L^\infty({\mathbb R}^N)}+
\sup_{s\in(0,1)}|||u_{ij}(s)|||_{\eta,\Psi;1}\\
&\le C\left(|||\mu|||_{\eta,\Psi;1}+j^{-1}\right)^{\frac{2}{\kappa_1}}
\end{split}
\end{equation}
for $i\ge 1$ and $j\ge j_*$. 
\vspace{3pt}
\newline
\underline{Step 2.} 
We complete the proof of Theorem~\ref{Theorem:1.2}. 
By \eqref{eq:4.18} we apply \cite{V}*{Theorem~7.1} to obtain the following:
\begin{itemize}
  \item  
  for any compact set $K\subset{\mathbb R}^N\times(0,1)$, 
  there exist $C>0$ and $\omega\in(0,1)$ such that
  \begin{equation*}
  |u_{ij}(x_1,t_1)-u_{ij}(x_2,t_2)|\le C\left(|x_1-x_2|^\omega+|t_1-t_2|^{\frac{\omega}{2}}\right)
  \end{equation*}
  for $(x_1,t_1)$, $(x_2,t_2)\in K$, $i\ge 1$, and $j\ge j_*$. 
\end{itemize} 
By the Arzel\`a-Ascoli Theorem and the diagonal argument 
we find a subsequence $\{u_{ij}'\}$ of $\{u_{ij}\}$ and a H\"older continuous function $u$ in ${\mathbb R}^N\times(0,1)$ 
such that 
$$
\lim_{i,j\to\infty}\|u_{ij}'-u\|_{L^\infty(K)}=0
$$
for any compact set $K$ of ${\mathbb R}^N\times(0,1)$. 
Then we observe from \eqref{eq:3.7}, \eqref{eq:3.21}, and \eqref{eq:4.22} that 
\begin{equation*}
\begin{split}
	\lim_{i, j\to\infty}\int_{B(z, 1)}u_{ij}(t)^p\,\dee x
	&=\int_{B(z, 1)}u(t)^p\,\dee x,\\
	\int_{B(z, 1)}u_{ij}(t)^p\,\dee x
	&\le \sup_{x\in\mathbb{R}^N}\left\{u_{ij}(x, t)^{p-1}\left(\log\left(e+u_{ij}(x, t)\right)\right)^{-\alpha}\right\}\int_{B(z,1)}\Psi(u_{ij}(t))\,\dee x\\
	& \le C\|u_{ij}(t)\|_{L^\infty(\mathbb{R}^N)}^{p-1}\left(\log\left(e+\|u_{ij}(t)\|_{L^\infty(\mathbb{R}^N)}\right)\right)^{-\alpha}\Psi\left(|||u_{ij}(t)|||_{\eta,\Psi;1}\right)\\
	& \le Ct^{-1}\left[\log\left(e+\frac{1}{t}\right)\right]^{-1-\alpha}\left(|||\mu|||_{\eta,\Psi;1}+1\right)^{\frac{2(p-1)}{\kappa_1}}\Psi\left(|||\mu|||_{\eta,\Psi;1}+1\right)
\end{split}	
\end{equation*}
for $z\in{\mathbb R}^N$, $t\in(0, 1)$, $i\ge 1$, and $j\ge j_*$.  
Therefore, by Definition~\ref{Definition:1.1}-(2), \eqref{eq:4.6}, and \eqref{eq:4.22} 
we apply the Lebesgue dominated convergence theorem to see that
$u$ is a solution to problem~\eqref{eq:P} in ${\mathbb R}^N\times(0,1)$ satisfying
\begin{equation*}
\sup_{s\in(0,1)} s^{\frac{1}{p-1}}\left[\log\left(e+\frac{1}{s}\right)\right]^{\frac{1}{p-1}}\|u(s)\|_{L^\infty({\mathbb R}^N)}+
\sup_{s\in(0,1)}|||u(s)|||_{\eta,\Psi;1}\le C|||\mu|||_{\eta,\Psi;1}^{\frac{2}{\kappa_1}}.
\end{equation*}
Thus Theorem~\ref{Theorem:1.2} follows.
$\Box$
\section{Proof of Theorem~\ref{Theorem:1.3}}
In this section we modify arguments in Section~4 to study sufficient conditions 
for the existence of solutions to problem~\eqref{eq:P} with $p>p_m$, 
and prove Theorem~\ref{Theorem:1.3}.
\vspace{5pt}
\newline
{\bf Proof of Theorem~\ref{Theorem:1.3}.}
Let $p>p_m$, $1<\beta<N(p-m)/2$, $T\in(0,\infty]$, and $\mu\in\mathcal{L}$. 
Let $\epsilon_6\in(0,1)$ be small enough, and assume \eqref{eq:1.6}.  
For any $i$, $j=1,2,\dots$, 
let $u_{ij}$ be a solution to problem~\eqref{eq:P} with initial data $\mu$ replaced by \eqref{eq:4.6}.
Then, for any $n\ge 1$, we find $j_*=j_*(n,\epsilon_6)$ such that 
\begin{equation}
\label{eq:5.2}
\begin{split}
|||\mu_{ij}|||_{\frac{N(p-m)}{2},\beta;T_n^{\theta}}
 & \le |||\mu|||_{\frac{N(p-m)}{2},\beta;T^{\theta}}+T_n^{\frac{2\theta}{p-m}}j^{-1}\\
 & \le |||\mu|||_{\frac{N(p-m)}{2},\beta;T^{\theta}}+
n^{\frac{1}{p-1}}j^{-1}\le 2\epsilon_6
\end{split}
\end{equation}
for $i\ge 1$ and $j\ge j_*$, where $T_n:=\min\{T,n\}$. 
Similarly to the proof of Theorem~\ref{Theorem:1.2}, 
by arguments in \cite{ADi} we find a unique classical solution $u_{ij}$ to problem~\eqref{eq:P} in ${\mathbb R}^N\times(0,T_{ij})$, 
with $u_{ij}$ satisfying \eqref{eq:4.7}, \eqref{eq:4.8}, and \eqref{eq:4.9}, 
where $T_{ij}$ is the maximal existence time of $u_{ij}$. 
\vspace{3pt}
\newline
\underline{Step 1.}
Let $n=1,2,\dots$, and fix it.
Let $\delta_1$, $\delta_2\in(0,1)$.
Set
\begin{align}
\label{eq:5.3}
T^1_{ij}:= & \,\sup\left\{t\in(0,T_{ij})\,:\, \sup_{s\in(0,t)}\,
s^{\frac{1}{p-1}}\|u_{ij}(s)\|_{L^\infty({\mathbb R}^N)}\le \delta_1\right\},\\
\label{eq:5.4}
T^2_{ij}:= & \,\sup\left\{t\in(0,T_{ij})\,:\,\sup_{s\in(0,t)}|||u_{ij}(s)|||_{\frac{N(p-m)}{2},\beta;T_n^{\theta}}\le\delta_2\right\}.
\end{align}
Under suitable choices of $\delta_1$ and $\delta_2$, 
taking $\epsilon_6>0$ small enough if necessary, 
we show that 
$$
T^*_{ij}:=\min\{T^1_{ij},T^2_{ij},T_n\}=T_n\quad\mbox{for $i\ge 1$ and $j\ge j_*$}. 
$$ 
In the proof of Theorem~\ref{Theorem:1.3}, 
the constants $C$ are independent of $n$, $i\ge 1$, and $j\ge j_*$.

Similarly to the argument in the proof of Theorem~\ref{Theorem:1.2}, 
we see that $T^*_{ij}>0$ for $i\ge 1$ and $j\ge j_*$.
By Lemma~\ref{Lemma:3.2} and \eqref{eq:5.4} 
we have
\begin{equation*}
\begin{split} 
 & \sup_{z\in{\mathbb R}^N}\int^t_0\int_{B(z,\sigma)}u_{ij}^{m+\beta-3}|\nabla u_{ij}|^2\,\dee x\,\dee s<\infty,\\
 & \sup_{s\in(0,t]}\sup_{z\in{\mathbb R}^N}\int_{B(z,\sigma)}u_{ij}(s)^\beta\,\dee x+\sup_{z\in{\mathbb R}^N}\int^t_0\int_{B(z,\sigma)} u_{ij}^{m+\beta-3}|\nabla u_{ij}|^2\,\dee x\,\dee s\\
 & \le C\sup_{z\in{\mathbb R}^N}\int_{B(z,\sigma)}\mu_{ij}^\beta\,\dee x
 +C\sigma^{-2}\sup_{z\in{\mathbb R}^N}\int^t_0\int_{B(z,\sigma)}u_{ij}^{m+\beta-1}\,\dee x\,\dee s\\
 & \qquad
 +C\delta_2^{p-m} 
 \sup_{z\in{\mathbb R}^N}\int^t_0\int_{B(z,\sigma)}u_{ij}^{m+\beta-3}|\nabla u_{ij}|^2\,\dee x\,\dee s,
\end{split}
\end{equation*}
for $t\in(0,T_{ij}^*)$, $\sigma>0$, $i\ge 1$, and $j\ge j_*$. 
Taking $\delta_2>0$ small enough if necessary, 
by \eqref{eq:5.3} we obtain
\begin{equation*}
\begin{split}
Y_{ij}(t)
 & \le C\sup_{z\in{\mathbb R}^N}\dashint_{B(z,\sigma)}\mu_{ij}^\beta\,\dee x
 +C\sigma^{-2}\int^t_0\|u_{ij}(s)\|_{L^\infty({\mathbb R}^N)}^{m-1}Y_{ij}(s)\,\dee s\\
  & \le C\sup_{z\in{\mathbb R}^N}\dashint_{B(z,\sigma)}\mu_{ij}^\beta\,\dee x
 +C\sigma^{-2}\int^t_0 s^{-\frac{m-1}{p-1}} Y_{ij}(s)\,\dee s
\end{split}
\end{equation*}
for $t\in(0,T_{ij}^*)$, $\sigma>0$, $i\ge 1$, and $j\ge j_*$,
where
$$
Y_{ij}(t):=\sup_{z\in{\mathbb R}^N}\dashint_{B(z,\sigma)}u_{ij}(t)^\beta\,\dee x.
$$
Then Gronwall's inequality  implies that 
\begin{align*}
 & Y_{ij}(t)
 \le C\exp\left(\sigma^{-2}\int_0^t s^{-\frac{m-1}{p-1}}\,\dee s\right)\sup_{z\in{\mathbb R}^N}\dashint_{B(z,\sigma)}\mu_{ij}^\beta\,\dee x\\
 & \qquad\quad
 \le C\exp\left(\sigma^{-2}t^{\frac{p-m}{p-1}}\right)\sup_{z\in{\mathbb R}^N}\dashint_{B(z,\sigma)}\mu_{ij}^\beta\,\dee x
 \le C\sup_{z\in{\mathbb R}^N}\,\dashint_{B(z,\sigma)}\mu_{ij}^\beta\,\dee x
\end{align*}
for $t\in(0,T_{ij}^*)$, $\sigma\ge t^{\theta}$, $i\ge 1$, and $j\ge j_*$. 
Therefore we deduce from \eqref{eq:5.2} that
\begin{equation}
\label{eq:5.5}
\sup_{s\in(0,t)}\sup_{z\in{\mathbb R}^N}\sup_{\sigma\in[s^{\theta},T_n^{\theta})}\sigma^{\frac{2}{p-m}}\left(\,\dashint_{B(z,\sigma)}u_{ij}(s)^\beta\,\dee x\right)^{\frac{1}{\beta}}\le C|||\mu_{ij}|||_{\frac{N(p-m)}{2},\beta;T_n^{\theta}}\le C\epsilon_6
\end{equation}
for $t\in(0,T_{ij}^*)$, $i\ge 1$, and $j\ge j_*$. 

On the other hand, taking $\delta_1>0$ small enough if necessary, 
by \eqref{eq:5.3} we have 
\begin{align*}
 & \left(\frac{s^{\theta}}{2}\right)^{-2}\|u_{ij}(s)\|^{m-1}_{L^\infty({\mathbb R}^N)}+\|u_{ij}(s)\|_{L^\infty({\mathbb R}^N)}^{p-1}\\
 & \le 4\delta_1^{m-1} s^{-\frac{p-m}{p-1}}s^{-\frac{m-1}{p-1}}+\delta_1^{p-1}s^{-1}
 \le C\delta_1^{m-1}s^{-1}+C\delta_1^{p-1}s^{-1}\le 2s^{-1}
\end{align*}
for $s\in(0,T_{ij}^*)$, $i\ge 1$, and $j\ge j_*$. 
Then, by Lemma~\ref{Lemma:3.3} with $\sigma=t^{\theta}/2$ and $r=\beta$,
taking $\epsilon_6>0$ small enough if necessary,  by \eqref{eq:5.2} and \eqref{eq:5.5} we have
\begin{equation}
\label{eq:5.6}
\begin{split}
 & t^{\frac{1}{p-1}}\|u_{ij}(t)\|_{L^\infty({\mathbb R}^N)}
 \le Ct^{\frac{1}{p-1}-\frac{N+2}{\kappa_\beta}}\left(\sup_{z\in{\mathbb R}^N}\int_{t/2}^t\int_{B(z,t^{\theta})} u_{ij}^\beta\,dy\,ds\right)^{\frac{2}{\kappa_\beta}}\\
 & \qquad
 \le Ct^{\frac{1}{p-1}-\frac{N}{\kappa_\beta}}t^{\frac{2N\theta}{\kappa_\beta}}\sup_{s\in(t/2,t)}
 \sup_{z\in{\mathbb R}^N}\left(\,\dashint_{B(z,t^{\theta})} u_{ij}(s)^\beta\,\dee y\right)^{\frac{2}{\kappa_\beta}}\\
 & \qquad
 \le Ct^{\frac{1}{p-1}-\frac{N}{\kappa_\beta}+\frac{2N\theta}{\kappa_\beta}-\frac{2\theta}{p-m}\frac{2\beta}{\kappa_\beta}}\sup_{s\in(t/2,t)}
 \sup_{z\in{\mathbb R}^N}\left[t^{\frac{2\theta}{p-m}}\left(\,\dashint_{B(z,t^{\theta})} u_{ij}(s)^\beta\,\dee y\right)^{\frac{1}{\beta}}\right]^{\frac{2\beta}{\kappa_\beta}}\\
 & \qquad
 \le C\sup_{s\in(0,t)}\sup_{z\in{\mathbb R}^N}\sup_{\sigma\in[s^{\theta},T_n^{\theta})}
 \left[\sigma^{\frac{2}{p-m}}\left(\,\dashint_{B(z,\sigma)} u_{ij}(s)^\beta\,\dee y\right)^{\frac{1}{\beta}}\right]^{\frac{2\beta}{\kappa_\beta}}\\
 & \qquad
 \le C|||\mu_{ij}|||_{\frac{N(p-m)}{2},\beta;T_n^{\theta}}^\frac{2\beta}{\kappa_\beta}\le C\epsilon_6^{\frac{2\beta}{\kappa_\beta}}\le \frac{\delta_1}{2}
\end{split}
\end{equation}
for $t\in(0,T_{ij}^*)$, $i\ge 1$, and $j\ge j_*$. 
Furthermore, we observe from \eqref{eq:5.2} and \eqref{eq:5.6} that
\begin{equation}
\label{eq:5.7}
\begin{split}
 & \sup_{s\in(0,t)}\sup_{z\in{\mathbb R}^N}\sup_{\sigma\in(0,s^{\theta})}
\sigma^{\frac{2}{p-m}}\left(\dashint_{B(z,\sigma)} u_{ij}(s)^\beta\,\dee y\right)^{\frac{1}{\beta}}\\
 & \le\sup_{s\in(0,t)}
s^{\frac{1}{p-1}}\|u_{ij}(s)\|_{L^\infty({\mathbb R}^N)}
\le C|||\mu_{ij}|||_{\frac{N(p-m)}{2},\beta;T_n^{\theta}}^\frac{2\beta}{\kappa_\beta}
 \le C\epsilon_6^{\frac{2\beta}{\kappa_\beta}}
\end{split}
\end{equation}
for $t\in(0,T_{ij}^*)$, $i\ge 1$, and $j\ge j_*$. 
By \eqref{eq:5.5} and \eqref{eq:5.7}, 
taking $\epsilon_6>0$ small enough if necessary, we obtain 
\begin{equation}
\label{eq:5.8}
\begin{split}
\sup_{s\in(0,t)}|||u_{ij}(s)|||_{\frac{N(p-m)}{2},\beta;T_n^{\theta}}
\le C\epsilon_6+C\epsilon_6^{\frac{2\beta}{\kappa_\beta}}\le\frac{\delta_2}{2}
\end{split}
\end{equation}
for $t\in(0,T^*_{ij})$, $i\ge 1$, and $j\ge j_*$. 
Therefore, thanks to \eqref{eq:5.6} and \eqref{eq:5.8}, 
by the definition of $T_{ij}^*$, 
for any $n=1,2,\dots$, 
we see that $T_{ij}^*=T_n$ for $i\ge 1$ and $j\ge j_*$. 
Furthermore, by \eqref{eq:5.2}, \eqref{eq:5.5}, \eqref{eq:5.6}, and  \eqref{eq:5.7} 
we have
\begin{equation}
\label{eq:5.9}
\begin{split}
 & \sup_{s\in(0,T_n)}\,s^{\frac{1}{p-1}}\|u_{ij}(s)\|_{L^\infty({\mathbb R}^N)}+\sup_{s\in(0,T_n)}|||u_{ij}(s)|||_{\frac{N(p-m)}{2},\beta;T_n^{\theta}}\\
 & \le C\left(|||\mu|||_{\frac{N(p-m)}{2},\beta;T_n^{\theta}}+n^{\frac{1}{p-1}}j^{-1}\right)^{\frac{2\beta}{\kappa_\beta}}
\end{split}
\end{equation}
for $i\ge 1$ and $j\ge j_*$. 
\vspace{3pt}
\newline
\underline{Step 2.} 
We complete the proof of Theorem~\ref{Theorem:1.3}. 
By \eqref{eq:5.9} we apply \cite{V}*{Theorem~7.1} to obtain the following:
\begin{itemize}
  \item  
  for any compact set $K\subset{\mathbb R}^N\times(0,T_n)$, 
  there exist $C>0$ and $\omega\in(0,1)$ such that
  $$
  |u_{ij}(x_1,t_1)-u_{ij}(x_2,t_2)|\le C\left(|x_1-x_2|^\omega+|t_1-t_2|^{\frac{\omega}{2}}\right)
  $$
  for $(x_1,t_1)$, $(x_2,t_2)\in K$, $i\ge 1$, and $j\ge j_*$.
\end{itemize} 
By the Arzel\`a-Ascoli Theorem and the diagonal argument 
we find a subsequence $\{u_{ij}'\}$ of $\{u_{ij}\}$ and a H\"older continuous function $u$ in ${\mathbb R}^N\times(0,T_n)$ 
such that 
$$
\lim_{i,j\to\infty}\|u_{ij}'-u\|_{L^\infty(K)}=0
$$
for any compact set $K$ of ${\mathbb R}^N\times(0,T_n)$. 
Then we observe from Jensen's inequality and \eqref{eq:5.9} that 
\begin{equation*}
\begin{split}
	\lim_{i, j\to\infty}\int_{B(z,T_n^\theta)}u_{ij}(t)^p\,\dee x
	&=\int_{B(z,T_n^\theta)}u(t)^p\,\dee x,\\
	\int_{B(z,T_n^\theta)}u_{ij}(t)^p\,\dee x
	&\le T_n^{N\theta}\left(\,\dashint_{B(z,T_n^\theta)}u_{ij}(t)^{\beta}\,\dee x\right)^{\frac{p}{\beta}}\\
	&\le T_n^{N\theta-\frac{p}{p-1}}|||u_{ij}(t)|||_{\frac{N(p-m)}{2},\beta;T_n^\theta}^p\\
	&\le C T_n^{N\theta-\frac{p}{p-1}}\left(|||\mu|||_{\frac{N(p-m)}{2},\beta;T_n^\theta}+n^\frac{1}{p-1}\right)^{\frac{2\beta p}{\kappa_\beta}}\quad\text{if}\quad 1<p\le\beta,\\
	\int_{B(z,T_n^\theta)}u_{ij}(t)^p\,\dee x
	&\le \|u_{ij}(t)\|_{L^\infty(\mathbb{R}^N)}^{p-\beta}\int_{B(z,  T_n^\theta)}u_{ij}(t)^{\beta}\,\dee x\\
	&\le C T_n^{N\theta-\frac{\beta}{p-1}}\|u_{ij}(t)\|_{L^\infty(\mathbb{R}^N)}^{p-\beta}|||u_{ij}(t)|||_{\frac{N(p-m)}{2},\beta;T_n^\theta}^\beta\\
	& \le CT_n^{N\theta-\frac{\beta}{p-1}}t^{-1+\frac{\beta-1}{p-1}}\left(|||\mu|||_{\frac{N(p-m)}{2},\beta;T_n^\theta}+n^\frac{1}{p-1}\right)^{\frac{2\beta p}{\kappa_\beta}}\quad\text{if}\quad 1<\beta<p,
\end{split}	
\end{equation*}
for $z\in{\mathbb R}^N$, $t\in(0, T_n)$, $i\ge 1$, and $j\ge j_*$.
Therefore, by Definition~\ref{Definition:1.1}-(2), \eqref{eq:4.6}, and \eqref{eq:5.9} 
we apply the Lebesgue dominated convergence theorem to see that 
$u$ is a solution to problem~\eqref{eq:P} in ${\mathbb R}^N\times(0,T_n)$ satisfying
\begin{equation}
\label{eq:5.10}
\sup_{s\in(0,T_n)}||u(s)|||_{\frac{N(p-m)}{2},\beta;T_n^{\theta}}+\sup_{s\in(0,T_n)}s^{\frac{1}{p-1}}\|u(s)\|_{L^\infty({\mathbb R}^N)}\le C|||\mu|||_{\frac{N(p-m)}{2},\beta;T_n^{\theta}}^{\frac{2\beta}{\kappa_\beta}}
\end{equation}
for $t\in(0,T_n)$. Since $n$ is arbitrary, \eqref{eq:5.10} holds with $T_n$ replaced by $T$. 
Thus $u$ is our desired solution to problem~\eqref{eq:P}, 
and Theorem~\ref{Theorem:1.3} follows.
$\Box$
\section{Proof of Corollary~\ref{Corollary:1.1}: Optimal singularity}
In this section, applying Theorems~\ref{Theorem:1.1}--\ref{Theorem:1.3}, we prove Corollary~\ref{Corollary:1.1}.
\vspace{5pt}

\noindent
{\bf Proof of Corollary~\ref{Corollary:1.1}.}
Let $p=p_m$, $\alpha\in(0,N/2)$, and 
$$
\mu(x)=|x|^{-N}\displaystyle{\biggr[\log\biggr(e+\frac{1}{|x|}\biggr)\biggr]^{-\frac{N}{2}-1}}
$$
for a.a.~$x\in{\mathbb R}^N$. 
Let $\Psi$ be as in Theorem~\ref{Theorem:1.2}. 
Then, by \eqref{eq:3.21} we have
\begin{equation}
\label{eq:6.1}
\begin{split}
\Psi(c\mu(x)) & =(c\mu(x))^{\frac{1}{2}}(c\mu(x))^{\frac{1}{2}}\left[\log\left(e+c\mu(x)\right)\right]^\alpha
\preceq c^{\frac{1}{2}}\mu(x)\left[\log\left(e+\mu(x)\right)\right]^\alpha\\
 & \preceq c^{\frac{1}{2}}|x|^{-N}\left[\log\left(e+\frac{1}{|x|}\right)\right]^{\alpha-\frac{N}{2}-1}
\end{split}
\end{equation}
for a.a.~$x\in{\mathbb R}^N$ and $c\in(0,1)$. This implies that
$$
\dashint_{B(z,\sigma)}\Psi(c\mu(x))\,\dee x\preceq\dashint_{B(0,3\sigma)}\Psi(c\mu(x))\,\dee x
\asymp
c^{\frac{1}{2}}\sigma^{-N}\biggr[\log\biggr(e+\frac{1}{\sigma}\biggr)\biggr]^{\alpha-\frac{N}{2}}
$$
for $z\in B(0,2\sigma)$, $\sigma\in(0,1)$, and $c\in(0,1)$. 
Thus we have
\begin{equation}
\label{eq:6.2}
\Psi^{-1}\left(\,\dashint_{B(z,\sigma)}\Psi(c\mu(x))\,\dee x\right)\preceq c^{\frac{1}{2}}\sigma^{-N}\biggr[\log\biggr(e+\frac{1}{\sigma}\biggr)\biggr]^{-\frac{N}{2}}
=c^{\frac{1}{2}}\eta(\sigma)^{-1}
\end{equation}
for $z\in B(0,2\sigma)$, $\sigma\in(0,1)$, and $c\in(0,1)$.
On the other hand, 
\begin{equation}
\label{eq:6.3}
\Psi^{-1}\left(\,\dashint_{B(z,\sigma)}\Psi(c\mu(x))\,\dee x\right)\preceq c\|\mu\|_{L^\infty(B(z,\sigma))}
\preceq c\sigma^{-N}\biggr[\log\biggr(e+\frac{1}{\sigma}\biggr)\biggr]^{-\frac{N}{2}}\asymp c\eta(\sigma)^{-1}
\end{equation}
for $z\in {\mathbb R}^N\setminus B(0,2\sigma)$, $\sigma\in(0,1)$, and $c\in(0,1)$. 
By \eqref{eq:6.2} and \eqref{eq:6.3} 
we apply Theorem~\ref{Theorem:1.2} with $T=1$ to see that 
problem~\eqref{eq:P} possesses a solution in ${\mathbb R}^N\times(0,1)$ if $c$ is small enough.
Thus assertion~(1) holds if $p=p_m$.

On the other hand, we have
$$
\dashint_{B(0,\sigma)}c\mu(y)\,\dee y\succeq c\sigma^{-N}\biggr[\log\biggr(e+\frac{1}{\sigma}\biggr)\biggr]^{-\frac{N}{2}}
$$
for $\sigma\in(0,1)$. 
Then it follows from Theorem~\ref{Theorem:1.1} that problem~\eqref{eq:P} possesses no local-in-time solution if $c$ is large enough. 
Thus assertion~(2) holds if $p=p_m$. 
Assertions in the case of $p>p_m$ follows from similar arguments to 
those in the case of $p=p_m$. Therefore the proof of Corollary~\ref{Corollary:1.1} is complete.
$\Box$


\end{document}